\documentclass[11pt]{article}
\usepackage{amssymb}
\usepackage{graphicx}
\usepackage{epsfig}
\usepackage{amsmath}
\usepackage{amsfonts}

\begin{document}
 \baselineskip 0.21in

\title{\textbf{A Theoretical Analysis of Sparse Recovery Stability of Dantzig Selector and LASSO  }}
\author{YUN-BIN ZHAO\thanks{School of Mathematics, University of Birmingham,
Edgbaston, Birmingham B15 2TT,  United Kingdom ({\tt
y.zhao.2@bham.ac.uk}).  The research of this author was  partially supported by
the Engineering and Physical Sciences Research Council (EPSRC) under
  grant \#EP/K00946X/1.} and DUAN LI\thanks{Department of Systems Engineering and Engineering Management,
The Chinese University of Hong Kong, Shatin, NT, Hong Kong ({\tt
dli@se.cuhk.edu.hk}). The research of this author was partially supported by the Research Grants Council of Hong Kong under grant
RGC \#2150841. }  }

\date{(1st version:  1 November 2016; 2rd version: April 25, 2017)}

\maketitle

  \noindent \textbf{Abstract.}  Dantzig selector (DS) and LASSO problems have attracted plenty of attention in statistical learning, sparse data recovery and mathematical optimization.   In this paper, we provide a theoretical analysis of the sparse recovery stability of these optimization problems in more general settings and from a new perspective.   We establish recovery error bounds for these optimization problems   under a mild assumption
 called weak range space property of a transposed design matrix. This assumption is less restrictive than the well  known sparse recovery conditions such as  restricted isometry property (RIP),  null space property (NSP)  or   mutual coherence. In fact, our analysis indicates that this assumption is tight and cannot be relaxed for the standard DS problems in order to maintain their sparse recovery stability.  As a result, a series of new stability results  for DS and LASSO have been established under various matrix properties, including the RIP with constant $\delta_{2k}< 1/\sqrt{2}$ and the (constant-free) standard NSP of order $k.$ We prove that these matrix properties can yield an identical recovery error bound for DS and LASSO with stability coefficients  being measured  by the so-called Robinson's constant, instead of the conventional RIP or NSP  constant. To our knowledge, this is the first time that the stability results with such a unified feature are established for DS and LASSO problems. Different from the standard analysis in this area of research, our analysis is carried out deterministically, and  the key analytic tools used in our analysis include the  error bound of linear systems due to Hoffman and Robinson and  polytope approximation of symmetric convex bodies due to Barvinok.  \\



\newpage

\section{Introduction}

 Let $\mathbb{R}^n_+ $ denote the nonnegative orthant of the  Euclidean space  $\mathbb{R}^n.$ The  set of $m\times n $ matrices will be denoted by $\mathbb{R}^{m\times n }.$ All vectors are column vectors  and the identity matrix of any
order will be denoted by $I, $ unless otherwise stated.
     We use $\phi: \mathbb{R}^q \to \mathbb{R}_+ $ to denote a general norm on $\mathbb{R}^q$ satisfying that $\phi(\textbf{e}_i) =1$ for all $i=1, \dots, q,$ where $\textbf{e}_i\in \mathbb{R}^q, i=1, \dots, q,$ are the standard basis of $\mathbb{R}^q,$ i.e., the column vectors of the  $q\times q$ identity matrix. In particular, we use $\|\cdot \|_p: \mathbb{R}^q \to \mathbb{R}_+ $ to denote the $\ell_p$-norm, i.e.,  $\|x\|_p = \left(\sum_{i=1}^q |x_i|^p \right)^{1/p}  $ for $x\in \mathbb{R}^q,$ where $ p\in [1, \infty].$
In particular, $\|x\|_1 =\sum_{i=1}^q |x_i| $  and    $\|x\|_\infty =\max_{1\leq i\leq q} |x_i|. $

  Given   matrices  $A\in \mathbb{R}^{m\times n} $ $(m<n)$  and $M\in \mathbb{R}^{m\times q} $ ($m\leq q$) with $\textrm{rank} (A)=
  \textrm{rank} (M)= m, $  we consider the following $\ell_1$-minimization problem:
\begin{equation} \label{L1} \min_ {x} \{\|x\|_1: ~ \phi(M^T(Ax-y)) \leq \tau \}, \end{equation} where $\tau \in \mathbb{R}_+ $ is  a given parameter and $y\in \mathbb{R}^m$ is a given vector. In  signal recovery scenarios, $A $
is often called a design or sensing matrix consisting of  a
set of known or learned dictionaries, and  $
y=A\widehat{x} +\theta $ is a measurement vector acquired for the signal
$\widehat{x} $ to  recover, and $\theta \in \mathbb{R}^m $ is the measurement
error, bounded as $\phi(M^T\theta) \leq \tau. $

  Problem (\ref{L1})  includes several important special cases.  In fact, when $\tau =0$,  the  problem   is reduced to the standard $\ell_1$-minimization, i.e., $\min \{\|x\|_1: Ax=y\},$  which is a signal recovery problem in noiseless situations (e.g., \cite{BDE09, E10, FR13}).     When $\tau >0$, $M=I$ and $\phi (\cdot)=\|\cdot \|_2, $  problem  (\ref{L1}) becomes a  nonlinear $\ell_1$-minimization   which has  been widely studied in the field of compressed sensing and signal processing (e.g., \cite{CDS98, D06,E10, EK12, FR13}). Taking $M=A$ in (\ref{L1}) leads to the problem
\begin{equation} \label{S-GDS} \min_ {x} \{\|x\|_1: ~ \phi(A^T(Ax-y)) \leq \tau \}. \end{equation}
When   $\phi(\cdot) = \|\cdot\|_\infty,   $  this problem  becomes the standard Dantzig selector (DS) introduced by Cand\`es and Tao  \cite{CT07, CT07b}.
In this paper,   problem (\ref{L1})  is still referred to as  the Dantzig  selector (DS) although it is more general than the standard one. Closely related to  (\ref{L1}) is the   problem
  \begin{equation} \label{GLASSO} \min_ {x} \{\phi(M^T(Ax-y)): ~   \|x\|_1 \leq \mu \}, \end{equation}
    where $\mu>0 $ is a given parameter. This problem also includes several important special cases.  When $M= I$ and $\phi(\cdot)=\|\cdot \|_2 $,  problem (\ref{GLASSO}) becomes the well known LASSO (least absolute shrinkage and selection operator) problem  introduced by Tibshirani \cite{T96}:
     \begin{equation} \label{LASSO} \min_ {x} \{\|Ax-y\|_2: ~   \|x\|_1 \leq \mu \}. \end{equation}
  In addition, setting $M=A$ in (\ref{GLASSO}) yields the model
  \begin{equation} \label{DANTZIG-LASSO} \min_ {x} \{\phi(A^T(Ax-y)): ~   \|x\|_1 \leq \mu \}, \end{equation} which clearly is related to  the  DS problem (\ref{S-GDS}).
In this paper,   problem (\ref{GLASSO}) is still called  the LASSO problem  despite the fact that it is more general than   (\ref{LASSO}).

        To possibly recover the  sparse data $\widehat{x} $ satisfying the bound $\phi(M^T(A\widehat{x}-y)) \leq \tau,  $   problem (\ref{L1})   seeks the $\ell_1$-minimizer $x$ that complies with the same bound $\phi(M^T(Ax-y)) \leq \tau, $   while  problem (\ref{GLASSO})   minimizes the error $\phi(M^T(Ax-y)) $ by assuming that the recovered data $x$   and the original data $\widehat{x}$ obey the same $\ell_1$-norm bound. Under the conditions that the optimal solution of a recovery problem is unique and $A$ satisfies certain strong properties,   the  current stability analysis checks whether the difference between the original data and the found solution of a recovery problem can stay under control.  Given a recovery problem, however, the optimal solution of the problem is not always unique from a mathematical point of view, and more importantly  the  matrix might satisfy a condition less restrictive than the existing assumptions. This motivates one to carry out the stability analysis  of a general recovery problem such as (\ref{L1}) and (\ref{GLASSO}) under a mild assumption, which may apply to a wider  range of situations.

    DS and LASSO are   popular in the statistics literature \cite{CT07, CT07b, EHT07,  FS07,  BRT09, BG11, HTW15}.  As pointed out in \cite{JRL09, BRT09, C12, MRY07, AR09},    DS and LASSO exhibit a similar behavior in  many situations, especially in a sparsity scenario.    Under the sparsity assumption and  certain matrix conditions such as  mutual coherence, restricted isometry property (RIP) or null space property (NSP),  the  standard $\ell_1$-minimization has been shown to be stable in sparse data recovery \cite{F04, D06, DET06, CT05, CRT06b, F14, CWX13, AS14, CZ14}.
    These results might be   valid for DS and LASSO  under suitable  assumptions. For instance, the mutual coherence introduced in  \cite{F04,  T04} in signal processing was shown to ensure the recovery stability of LASSO  in \cite{MB06, ZY06}.   Under the RIP assumption,  Cand\`es and Tao \cite{CT07} have shown that if the true vector $\widehat{x}$ is sufficiently sparse  then   $\widehat{x}$ can be estimated reliably via DS based on the noise data observation $y$.  Cai, Xu and Zhang \cite {CXZ09} have  shown certain improved stability results for  DS and   LASSO by slightly weakening the condition in \cite {CT07}.  It is worth mentioning that the stability of LASSO can be guaranteed under the  restricted eigenvalue condition (REC) introduced by Bickel, Ritov and Tsybakov \cite{BRT09}. This condition holds  with high probability for    sub-Gaussian design matrices \cite{RZ13}.  Other stability conditions    have also been examined in the literature such as the compatibility condition \cite{GB09}, $H_{s,q}$ condition \cite{JN11}, and certain variants of RIP, NSP or REC \cite{C12,J15,KZL15, MR15}. A good summary of  stability results of LASSO can be found in \cite{BG11, HTW15}.

    A typical feature of the current stability theory  for DS and LASSO is that  the  coefficients in  error bounds are usually measured by  the RIP,  REC or other individual matrix constants. In the current framework, deferent matrix assumptions require distinct analysis and yield   different stability constants determined by the assumed matrix constants, such as the RIP constant, which are usually hard to certify (see, e.g., \cite{BDMS13, TP14}).
    The purpose of this  paper is to develop stability results for  general problems (\ref{L1}) and (\ref{GLASSO}) under a   constant-free matrix condition which is a very mild assumption.  In fact, it turns out that this assumption is both necessary and sufficient for the standard DS to be stable in sparse data recovery.   Our stability results and analysis   are developed and carried out under less restrictive assumptions than the existing ones and for more general optimization problems. More specifically,
our work is carried out differently in three aspects.

 (i) The results in this paper are established in a fairly general setting. Our analysis is based on the fundamental Karush-Kuhn-Tucker (KKT) optimality conditions \cite{KT51, K39} which capture  the deep  property of the optimal solution  of a convex optimization problem. Thus KKT conditions naturally lead to  the so-called \emph{weak  range space property (weak RSP) of order $k$ of $A^T$} (see Definitions 2.6 for details) which turns out to be a very  mild assumption guaranteeing the  stability of a broad range of sparse optimization problems including the standard   DS and LASSO. For the standard DS, we will show that this assumption is tight and cannot be  relaxed in order to guarantee its sparse recovery stability.

 (ii) The weak RSP of order $k$ of $A^T$ is a constant-free matrix property. A unique feature of the stability results developed in this paper is that the  recovery error bounds are  measured through so-called Robinson's constant \cite{R73} which depends on the given problem data.
Many known matrix conditions, such as the RIP \cite{CT05}, NSP \cite{CDD09}, mutual coherence \cite{E10}, REC \cite{BRT09},   and the range space property (RSP) of $A^T$ \cite {Z13, ZX14},     imply the weak RSP introduced in this paper. We show that the existing conditions imply an identical recovery error bound measured  with
Robinson's constants. Thus our results can be viewed as certain unified stability theorems for DS and LASSO problems.

(iii) The classic error bound of linear systems  (developed   by Hoffman \cite{H52}, and refined later by
Robinson \cite{R73} and Mangasarian \cite{MS87})  is a major analytic tool used in this paper. It provides a useful means for the study of stability of DS and LASSO problems. The stability of linear DS problems  can be established directly by such an error bound of linear systems.  Combined with the  polytope approximation of symmetric convex bodies developed by Barvinok  \cite{B14} (see also the earlier work by Dudley \cite{D74},   Bronshtein and  Ivanov \cite{BI75}, and Pisier \cite{P89}), the Hoffman's error bound can also be  used to establish the stability of nonlinear convex optimization problems, including the nonlinear DS and LASSO.

This paper is organized as follows. In Section 2, we provide basic facts, definitions, and initial results that will be used in the remainder of the paper.
The  stability of the linear DS under the weak RSP is proved in Section 3. Based on the polytope approximation of the unit ball,   we show in Sections 4 and 5  that the nonlinear problems (\ref{L1}) and (\ref{GLASSO}) are also stable  under the weak RSP.

  \section{Preliminaries}

 \subsection{Notations} In addition to the notations introduced at the beginning of Section 1, the following notations  will also be used throughout the   paper.     $x \in \mathbb{R}^n $ is called a $k$-sparse vector if it
admits at most $k$ nonzero entries,  i.e., $\| x\|_0 \leq k,$ where $\|\cdot\|_0$ counts the number of nonzero components of $x.$  For $x\in \mathbb{R}^n,$  let $|x|$,  $(x)^+$ and $(x)^-$
  be the vectors   with components $|x|_i
:= |x_i| $,   $ [(x)^+]_i := \max \{x_i, 0\}  $  and $ [(x)^-]_i := \min \{x_i, 0\}  $ for $ i=1, \dots, n,$ respectively.  For $x, y \in \mathbb{R}^n$,  the inequality $x\leq y $ means $x_i\leq y_i$ for
all $i=1, \dots, n.$      Given a set $S\subseteq
\{1, \dots, n\},$        we use $\overline{S} =\{1, \dots , n\}\backslash S$
    to denote
the complement of $S$ with respect to $\{1, \dots, n\}.$
We use  $x_S$ to denote the subvector of $x \in \mathbb{R}^n$ by deleting
the components $x_i$ with $i \in \overline{S},   $ and we use $Q_S$ to denote the submatrix of $Q\in \mathbb{R}^{m\times n} $  by removing
those columns of $Q$ with column index $i \in \overline{S}.  $        $Q^T$
denotes the transpose of the matrix $Q, $ and  ${\cal R} (Q^T)=  \{ Q^Tu:  u \in R^m\} $ denotes the range space of
$Q^T. $  For any given norms $\phi' : \mathbb{R}^n \to \mathbb{R}_+$ and $ \phi'' : \mathbb{R}^m\to \mathbb{R}_+$, the induced
  norm  $\|Q\|_{\phi' \to \phi'' }$ of the matrix $Q\in \mathbb{R}^{m\times n} $ is defined as $$\|Q\|_{\phi' \to \phi'' } = \max_{\phi'(x)\leq
1} \phi''(Qx).$$ In particular,  $\|Q\|_{p\to q} = \max_{\|x\|_p\leq
1} \|Qx\|_q $ where $p,q\geq 1.$  Given two sets $\Omega_1, \Omega_2
\subseteq \mathbb{R}^m,  $
we use $d^{\cal H} (\Omega_1, \Omega_2)$ to denote the Hausdorff distance of $(\Omega_1, \Omega_2)$, i.e.,
\begin{equation} \label{HM} d^{\cal H} (\Omega_1, \Omega_2) =
\max \left\{\sup_{u\in \Omega_1} \inf_{z\in \Omega_2} \|u-z\|_2, ~
\sup_{z\in \Omega_2} \inf_{u\in \Omega_1} \|u-z \|_2\right \}.
\end{equation}  For a closed convex set $\Omega \subseteq
\mathbb{R}^n$,  let $\Pi_ {\Omega} (x)$ denote the orthogonal projection of $x \in \mathbb{R}^n$ into $ \Omega, $ i.e., $
\Pi_ {\Omega} (x) := \textrm{arg}\min_{u} \{\| x-u\|_2:  ~ u\in \Omega\} . $

\subsection{Robinson's Constant and Hoffman's Lemma}
 Let $M'  \in
\mathbb{R}^{p_1\times q}$ and $ M'' \in \mathbb{R}^{p_2\times q} $ be
two   matrices. Let $\|\cdot\|_{\alpha_1}$ and $ \|\cdot\|_{\alpha_2}, $ where $\alpha_1, \alpha_2\geq 1,$ be the $\ell_{\alpha_1}$- and $\ell_ {\alpha_2} $-norms on
$\mathbb{R}^q$ and $ \mathbb{R}^{p_1+p_2}, $ respectively.  Define
    \begin{equation}
\label{HR-mu} \mu_{\alpha_1, \alpha_2} (M', M''):= \max_{(d',d'')\in {\cal F}, \|(d',d'')\|_{\alpha_2}
\leq 1}
 \min_{u\in \mathbb{R}^q} \{\|u\|_{\alpha_1}:  ~ M' u\leq d',  ~ M'' u =d''\},
 \end{equation} where ${\cal F}$ is the set defined by  $${\cal F} = \{(d',d'') \in \mathbb{R}^{p_1+p_2}:\textrm{ for some } u \in \mathbb{R}^q \textrm{ such that }  M' u \leq d' \textrm{ and }  M''  u =d''\}.$$
  Note that $\mu_{\alpha_1, \alpha_2} (M', M'')$ is a finite number (see Robinson \cite{R73}). Let $M^{1} \in \mathbb{R}^{m_1 \times \ell}$ and $ M^2 \in \mathbb{R}^{m_2
\times \ell }   $ be a  pair of matrices and   $ S$ be a subset of $    \{1, \dots, m_1\}. $  We consider the following two matrices:  $$  \left[
\begin{array}{cc}
 I_S & 0 \\
   -I & 0 \\
    \end{array}
    \right] \in \mathbb{R}^{(|S|+m_1)\times (m_1+m_2)},  ~  \left[
          \begin{array}{c}
     M^1 \\
      M^2\\
     \end{array}
     \right]^T \in \mathbb{R}^{\ell \times ( m_1+m_2 )},$$ where  $I_S$ is the submatrix extracted from the
       identity matrix $I \in \mathbb{R}^{m_1\times m_1} $ by deleting the rows corresponding to indices not in $S.$ Substituting the above pair into (\ref{HR-mu}) and taking the maximum over all possible subsets $S$ leads to
     \begin{equation} \label{HR-constant} \sigma_{\alpha_1, \alpha_2}
(M^1, M^2):= \max_{S\subseteq \{1, \dots , m_1\} } \mu_{\alpha_1, \alpha_2}
\left(\left[
                \begin{array}{cc}
     I_S & 0 \\
       -I & 0 \\
        \end{array}
         \right], \left[
          \begin{array}{c}
               M^1 \\
            M^2\\
           \end{array}
         \right]^T \right),
   \end{equation}
   which is a constant introduced by Robinson \cite{R73}. We call (\ref{HR-constant}) the Robinson's constant in this paper. Using this constant with $(\alpha_1,\alpha_2) =(\infty, 2)  $, Robinson \cite{R73}
  proved that the classic Hoffman's Lemma \cite{H52} concerning linear systems can be restated as follows.

 \vskip 0.07in

   \textbf{Lemma 2.1.} \cite{H52, R73} \emph{Let $M^1\in \mathbb{R}^{m_1\times \ell} $ and $ M^2 \in \mathbb{R}^{m_2\times \ell}$ be given matrices and ${\cal L}
= \{u\in \mathbb{R}^\ell : M^1u \leq d^1, M^2 u= d^2\} $ where $d^1\in \mathbb{R}^{m_1}$ and $ d^2\in \mathbb{R}^{m_2}.$
 Then for any   $x\in\mathbb{R}^\ell, $  there is a point $x^* \in  {\cal L}$
 such that }
  $$ \|x-x^*\|_2 \leq  \sigma_{\infty, 2}(M^1, M^2)
   \left\| \left[
    \begin{array}{c}
     (M^1x-d^1)^+ \\
      M^2x-d^2  \\
      \end{array}
      \right] \right\|_1.  $$

 \subsection{Polytope approximation of the unit ball}

Given a norm $ \phi(\cdot) $ on $\mathbb{R}^q,$   let $ \phi^*(\cdot) $ be the dual norm of $\phi,$ i.e.,  $ \phi^*(u) = \max_{\phi(x)\leq 1} u^Tx .$  From the  definition, we see that    $  x^T u \leq \phi(x) \phi^*(u) $  for any   $ x, u \in \mathbb{R}^q   .$
    In particular, the dual norm of  the $\ell_p$-norm, $ p \in [1, \infty] , $  is  the $\ell_\beta $-norm with $\beta \in [1, \infty] , $ where $p$ and $\beta$ satisfy that $1/p+ 1/\beta =1.  $ For instance, the $\ell_1$-norm   and $\ell_\infty $-norm   are dual to each other. In this paper, we restrict our attention to the norm $\phi$ with $\phi(\textbf{e}_i ) =1$ and $\phi^*(\textbf{e}_i ) =1$ for all $ i =1, \dots, q. $  Clearly, any $\ell_p$-norm  satisfies this property. Let $$ \mathfrak{B}^{\phi} = \{x \in \mathbb{R}^q : ~ \phi(x) \leq 1\}$$  be the unit ball in $ \mathbb{R}^q$  defined by  the norm $\phi   .$  It is evident that
  $$ \mathfrak{B}^{\phi}  =\bigcap_{\phi^*(a) =1} \{x \in \mathbb{R}^q :  a^T x \leq 1\}. $$  This means that   $\mathfrak{B}^{\phi}$ is the intersection of  all half spaces in the form $\{x \in \mathbb{R}^q:  a^T x \leq 1\}  $ where $a\in \mathbb{R}^q $ and $\phi^*(a) =1.$   Any finite number of  the vectors $a^i\in \mathbb{R}^q $ with $ \phi^* (a^i) =1$, $ i=1, \dots, k ,$  yield the following polytope approximation   of $ \mathfrak{B}^{\phi}:$  $$\mathfrak{B}^{\phi} \subseteq \bigcap_{1\leq i\leq k} \{x \in \mathbb{R}^q :  (a^i)^T x \leq 1\} . $$
     Note that $1=\phi^*(a^i)=\sup_{\phi(u)\leq 1} (a^i)^T u = (a^i)^T u^*$ for some $u^*$ with $\phi(u^*) =1.$     Thus every half space $ \{x \in \mathbb{R}^q:  (a^i)^T x \leq 1\}$ with $\phi^* (a^i) =1$ must be a support half space of $ \mathfrak{B}^{\phi}  $ in the sense that it contains $\mathfrak{B}^{\phi} $ and the plane $\{x: (a^i)^T x = 1\} $ touches $\mathfrak{B}^{\phi} $ at a point on its boundary.  Conversely, any support half space of $ \mathfrak{B}^{\phi} $ can be represented this way, i.e., $\{x: (a^i)^T x \leq 1 \} $ with some $a^i$ satisfying $\phi^* (a^i) =1. $   Note that   $ \mathfrak{B}^{\phi} $ is a symmetric convex body (i.e., if $x$ is in the set, so is $-x$), to which there is a polytope approximation \cite{D74, B14}, as claimed by the following theorem.

 \vskip 0.07in

\textbf{Theorem 2.2.} (Barvinok \cite{B14}) \emph{For any constant $\chi > \frac{e}{4 \sqrt{2}} \approx0.48, $ there exists an $\epsilon _0=\epsilon_0(\chi ) $ such that for any $0< \epsilon < \epsilon_0 $ and for any symmetric convex body $  B $ in $ \mathbb{R}^q,$ there is a symmetric polytope in $ \mathbb{R}^q $, denoted by $P_\epsilon,$  with $N$ vertices such that    $N\leq \left(\frac{\chi}{\sqrt{\epsilon}} \ln \frac{1}{\epsilon} \right)^q $    and  $P_\epsilon \subseteq  B  \subseteq (1+\epsilon) P_\epsilon. $}

 \vskip 0.08in

The above theorem indicates that   for every sufficiently  small $\epsilon >0$ there exists a   polytope $P_\epsilon$ satisfying
  $P_\epsilon \subseteq  \mathfrak{B}^\phi \subseteq (1+\epsilon) P_\epsilon. $ Thus $(1+\epsilon) P_\epsilon  $  is an outer approximation of $ \mathfrak{B}^\phi. $ To get a tighter approximation of  $\mathfrak{B}^\phi, $ we  compress the polytope $(1+\epsilon) P_\epsilon$   by shifting all affine planes, expanded by the faces of $(1+\epsilon) P_\epsilon,$  toward
$\mathfrak{B}^\phi $ until they touch  $\mathfrak{B}^\phi $ on its boundary. By such a compression, the resulting polytope   denoted by $ \widehat{P}_\epsilon  $  is then formed by a finite number of support half spaces of $\mathfrak{B}^\phi.$  Therefore, there exists a set  of vectors $a^i$ for $i=1 , \dots, J $ with   $\phi^*(a^i)=1  $    such that
$$  \widehat{P}_\epsilon =  \bigcap_{i=1}^{J}   \left\{x: (a^i)^T x  \leq 1\right\}.$$
 We now   add the   $2m$
half spaces $(\pm \textbf{e}_i)^T z \leq 1, ~i=1, \dots, q, $ to $\widehat{P}_\epsilon , $   yielding the following polytope:
\begin{equation}  \label {PK-C}  \overline{P}_\epsilon: = \widehat{P}_\epsilon \cap  \left\{z\in \mathbb{R}^q:   ~ \textbf{e}_i^T z\leq 1,
~ -\textbf{e}_i^Tz \leq  1, ~ i=1, \dots, q\right\}.
 \end{equation}
Clearly, we have the  relation:
$ P_\epsilon \subseteq \mathfrak{B}^\phi \subseteq  \overline{P}_\epsilon  \subseteq  \widehat{P}_\epsilon \subseteq (1+\epsilon) P_\epsilon . $
Throughout the   paper, we let $\epsilon_k \in (0, \epsilon_0), $ where $\epsilon_0$ is the constant specified in Theorem 2.2,  be a positive and strictly decreasing sequence
  satisfying   $\epsilon_k \to 0 $ as $k \to \infty.$       We consider the sequence of
polytopes $ \{ {\cal Q}_ {\epsilon_j}\}_{j\geq 1 } ,  $  where
 \begin{equation} \label {PJ-01}  {\cal Q}_{\epsilon_j}  = \bigcap_{k=1 }^ j \overline{P}_{\epsilon_k}  .
\end{equation} Then for every $\epsilon_j $,  $ {\cal Q}_{\epsilon_j} $ satisfies that
\begin{equation} \label{contain}  P_{\epsilon_j}  \subseteq \mathfrak{B}^\phi \subseteq   {\cal Q}_{\epsilon_j} \subseteq \overline{P}_{\epsilon_j} \subseteq (1+\epsilon_j) P_{\epsilon_j} . \end{equation}
Note that $ {\cal Q}_{\epsilon_j}$ is   a polytope formed by a
finite number (say, $K$) of half spaces, denoted by $(a^i)^T z \leq 1$ where $\phi^* (a^i ) =1 $ for $i=1, \dots, K. $
 We use  $ \Gamma_{ {\cal Q}_{\epsilon_j}}:
=\left [a^1, \dots , a^K\right]$ to denote the matrix with the vectors $a^i$'s as its
columns,  and we use   $\mathcal{Y} (\Gamma_{ {\cal Q}_{\epsilon_j}})
=\left\{a^1,  \dots, a^K\right\}  $   to denote the set of columns of
$\Gamma_{ {\cal Q}_{\epsilon_j}}. $    Then $ {\cal Q}_{\epsilon_j}$ can be expressed as \begin{equation} \label{QQQ}  {\cal Q}_{\epsilon_j}  = \left\{z\in \mathbb{R}^q: ~ (a^i)^T z\leq 1, ~ i=1,\dots, K \right\} =
\left\{z\in \mathbb{R}^q: (\Gamma_{ {\cal Q}_{\epsilon_j}})^T z\leq \widehat{\mathbf{e}} \right\} , \end{equation}  where $ \widehat{\mathbf{e}} $ denotes
the vector of ones in $\mathbb{R}^K. $ The following lemma is useful in our late analysis.

\vskip 0.08in

\textbf{Lemma 2.3.}  \emph{ For any $j\geq 1,$  let $ {\cal Q}_{\epsilon_j}$ be constructed as (\ref{PJ-01}).  Then for any   $a^* $ on the
 unit sphere $   \{x\in \mathbb{R}^q: \phi(x) =1\},$ there exists a   vector   $a^i \in
  {\cal Y}(\Gamma_{ {\cal Q }_{\epsilon_j}}) $  such that
 $   (a^*)^T a^i   \geq  \frac{1}{1+\epsilon_j}    .$}

\vskip 0.08in

\emph{Proof.} Let $a^*$ be any point on the unit sphere, i.e.,   $\phi(a^*)=1.$ Note that $ {\cal Q }_{\epsilon_j}$   satisfies (\ref{contain}).   The straight line starting from the  origin and passing through the point $a^*$ on the surface of $\mathfrak{B}^\phi $ will shoot a point $z'$ on the boundary of $ {\cal Q }_{\epsilon_j} $ and a point $z''$ on the boundary of $(1+\epsilon_j) P_{\epsilon_j }. $  From   (\ref{contain}), we see  that $z''= (1+\epsilon'') a^*$ for some number $\epsilon'' \leq \epsilon_j.$ Note that $z'$ is situated between $a^*$ and $ z''$. This implies that $z'= (1+\epsilon')a^*$ for some $ \epsilon'\leq \epsilon'' .$   Since $z'$ is on the boundary of $ {\cal Q}_{\epsilon_j}, $   it must be on a face of this polytope and hence there exists a vector $a^i \in  {\cal Y}(\Gamma_{ {\cal Q }_{\epsilon_j}})$ such that $ (a^i)^T z'=1.$ Note that $z'=(1+\epsilon') a^*$ where $ \epsilon'\leq \epsilon'' \leq \epsilon_j  . $ Thus
$ 1= (a^i)^T z' = (1+\epsilon') (a^i)^T a^* $ which implies that
$ (a^i)^T a^* = \frac{1}{1+\epsilon'}  \geq  \frac{1}{1+\epsilon_j}.$   ~~ $\Box$

\subsection{Stability and weak RSP condition}

Let us first give the definition of the stability of a sparse optimization problem. For a given vector $x\in \mathbb{R}^n,$  we recall that  the   best $k$-term
approximation of   $x  $   is defined as follows:
$$\sigma_k(x)_1 := \inf_{u} \{\|x-u\|_1: \|u\|_0\leq k\},$$
where $k$ is a given integer number   and $\|u\|_0 $ counts the number of nonzero entries of $u \in  \mathbb{R}^n.$   For   problem (\ref{L1}), the stability can be described as follows.

\vskip 0.05in

\textbf{Definition 2.4.}    Let $\widehat{x}$ be the original data obeying the constraint  of (\ref{L1}). Problem (\ref{L1})  is said to be \emph{stable} in   data recovery if there is an optimal solution $x^*$   of
 (\ref{L1}) approximating $\widehat{x}$ with error
\begin{equation} \label{RS} \|\widehat{x}-x^*\|_2 \leq  C'  \sigma_k(\widehat{x})_1 +C''
\tau, \end{equation}
where $C'$ and $C''$   are two  constants depending only on the problem data  $(M,A, y, \tau).$

\vskip 0.05in

For the LASSO problem  (\ref{GLASSO}),  we introduce the following definition.

\vskip 0.05in

\textbf{Definition 2.5.}   Let $\widehat{x}$ be the original data  obeying the constraint  of (\ref{GLASSO}).  Problem (\ref{GLASSO})  is said to be \emph{stable} in   data recovery if there is an optimal solution $x^*$   of
(\ref{GLASSO}) approximating  $\widehat{x}$ with error
\begin{equation} \label{RS1} \|\widehat{x}-x^*\|_2 \leq  C_1 \sigma_k(\widehat{x})_1 + C_2  (\mu-\|\widehat{x}\|_1) +  C_3 \mu \phi(M^T (A\widehat{x}-y)), \end{equation}
where   $  C_1, C_2 $ and $ C_3$   are  constants depending only on  the problem data $(M,A, y, \mu).$

\vskip 0.05in

We now introduce the weak range space property of $A^T.$   By the KKT optimality condition,   a $k$-sparse vector $\widehat{x}$ is an optimal
solution to the standard $\ell_1$-minimization problem $\min\{ \|x\|_1: Ax=y:= A\widehat{x}\}$ if and only if there is a vector $\zeta  \in
{\cal R}( A^T)
  $ satisfying that  $ \zeta_i =1 $  for $\widehat{x}_i>0,$   $\zeta_i= -
1  $ for $\widehat{x}_i<0 ,$  and $ |\zeta_i|\leq  1 $ for $\widehat{x}_i =0.$
 This property of $A^T$ depends on the individual vector $\widehat{x}.$  To ensure every $k$-sparse vector can be exactly recovered by $\ell_1$-minimization, it is necessary to
strengthen this property so that it is independent of any individual  $k$-sparse vector. This naturally leads to the next definition.

\vskip 0.07in

\textbf{Definition 2.6.} (Weak RSP of order $k$ of $A^T$)  Let $A $ be an $m\times n $ matrix with  $m<n. $   We say that $A^T$  admits
 the weak range space property of order $k$ if for any
disjoint subsets $S_1, S_2 \subseteq \{1,\dots, n\}$ with
$|S_1|+|S_2|\leq
 k$, there exists a vector $\zeta  \in {\mathcal R}(A^T) $ obeying
 \begin{equation}\label{RSP} \zeta_i =1  \textrm{ for   }  i \in S_1, ~  \zeta_i =-1
  \textrm{ for  }  i\in
 S_2, \textrm{  and   } |\zeta_i| \leq 1 \textrm{  for   } i\notin S_1\cup S_2. \end{equation}

\vskip 0.07in

  Slightly strengthening the condition ``$|\zeta_i| \leq 1 \textrm{  for   } i\notin S_1\cup S_2 $" to  ``$|\zeta_i| <  1 \textrm{  for   } i\notin S_1\cup S_2, $"    the above concept becomes the  RSP  of order $k$ of $A^T$  introduced  in \cite{Z13} (some earlier related works can be found in \cite{P07, GSH11}). It was shown in \cite{Z13} that the RSP of order $k$ of $A^T$ is a necessary and sufficient condition  for the  recovery of every $ k$-sparse signal  via the standard $\ell_1$-minimization.  Many  sparse recovery conditions must  imply the weak RSP. To see this, let us first
  recall  a few existing  matrix properties.  $A\in \mathbb{R}^{m\times n} $ is said to satisfy the restricted isometry property (RIP)
of order $k$ if there exists a constant $\delta_{k} \in (0,1)$ such that
$ (1-\delta_{k}) \|x\|^2_2 \leq \|Ax\|^2_2
\leq (1+\delta_{k})\|x\|^2_2$ for any  $k$-sparse vector
$x\in \mathbb{R}^n  $ (see \cite{CT05}).  $A$ satisfies the null space
property (NSP) of order $k$   if
$\|\zeta_S\|_1 < \|\zeta_{\overline{S}}\|_1
 $ holds for any $\zeta \not= 0$ in the null space of $A $ and any $S\subseteq \{1, \dots, n\}$ with $|S| \leq
k$ (see \cite{CDD09, FR13}).
Strengthening the NSP concept   leads to the following
 stable or robust NSP of order $k$ (see \cite{CDD09, FR13}):
(i)  $A\in \mathbb{R}^{m\times n} $  satisfies the \emph{stable null space property of order $k$}
if these is a constant $\rho \in (0, 1)$ such that  $\|\zeta_S\|_1 \leq \rho
\|\zeta_{\overline{S}}\|_1
 $   for any  $\zeta \not= 0$ in the null space of $A $ and any $S\subseteq \{1, \dots, n\}$ with $|S| \leq
k; $  (ii)  $A$ is said to admit the \emph{robust null space property of order $k$} if there are constants $\rho' \in (0, 1)$ and $\rho'' >0$ such that $\|\zeta_S\|_1 \leq \rho'
\|\zeta_{\overline{S}}\|_1 +\rho'' \|A \zeta\|
 $  for any  $\zeta \not= 0$ in the null space of $A $ and any $S\subseteq \{1, \dots, n\}$ with $|S| \leq
k. $
It is well known  that the
  mutual coherence condition $\mu_1(K)+\mu_1(k-1) <1$ introduced in \cite{T04, DET06, E10}  implies the RIP and NSP (see Theorem 5.15 in
\cite{FR13} and Lemma 1.5 in \cite{EK12}). The NSP is strictly
weaker than the RIP (e.g.  \cite{F14, CCW15}).
Note that  NSP of order $k$ is also  a necessary and sufficient condition for the  recovery of every $k$-sparse signal  (see \cite{FR13}),   the NSP of order $k$ is equivalent to the RSP of order $k, $  and hence each of them implies the weak RSP of order $k$ of $A^T.$
From the above discussion, we see that the RIP of $A$ is a strictly stronger than the RSP of $A^T.$  Thus many existing sparse recovery conditions imply the weak RSP of order $k$ of $A^T$ which, in fact, is a necessary condition for many recovery problems  to be stable as shown by the next theorem.

 \vskip 0.07in

\textbf{Theorem 2.7.} \emph{Let $\varphi: \mathbb{R}^q \to \mathbb{R}_+$ be a finite convex function on $\mathbb{R}^q$  satisfying  $\varphi(0) =0$ and $\varphi(u)>0$ for $u\not=0.$ Let $A\in \mathbb{R}^{m\times n}$ and $M\in \mathbb{R}^{m\times q}$, where $m< n $ and $ m\leq q, $ be  full-rank matrices.
Suppose that for any given $\varepsilon\geq 0$ and any given $y\in \mathbb{R}^m$,   the vector $x\in \mathbb{R}^n$ satisfying
$\varphi(M^T(Ax-y)) \leq \varepsilon $ can be approximated by an optimal solution $x^*$ of the problem \begin{equation}  \label{GP} \min_z \{\| z\|_1:  \varphi(M^T(Az-y)) \leq \varepsilon\} \end{equation} with error
 \begin{equation}  \label{GERROR} \|x-x^*\|_2 \leq  C' \sigma_k(x)_1 + C'' \gamma(\varepsilon), \end{equation}
where  $C'$ and $ C''$ are constants depending on the problem data $(M, A, y, \varepsilon),$  and $\gamma(\cdot)$ is a certain continuous function satisfying $\gamma(0) =0 ,  \gamma(\varepsilon) >0\textrm{  for } \varepsilon>0 , $ and $ C''\gamma(\varepsilon) \to 0\textrm{ as }\varepsilon \to 0.$       Then $A^T$ must satisfy the weak RSP of order $k.$}

\vskip 0.07in

  \emph{Proof.} Let $S_1, S_2 \subseteq \{1, \dots, n\}$  be  two  arbitrary and
disjoint sets with $|S_1|+|S_2| \leq k. $
We   prove that there is
a vector $\zeta \in {\cal R}(A^T) $ satisfying (\ref{RSP}). Let $\widehat{x}$ be a $k$-sparse
vector in $\mathbb{R}^n$ with
\begin{equation} \label{S1S2}  \{i: ~ \widehat{x}_i>0\}= S_1,  ~ \{i: ~ \widehat{x}_i<0\}= S_2. \end{equation}    Consider the  small parameter $ \varepsilon $     such that
\begin{equation} \label{choice}  C''\gamma(\varepsilon) <\min_{\widehat{x}_i\not =0} |\widehat{x}_i|.  \end{equation}  Let the measurements $ y \approx A\widehat{x} $ be accurate enough such that $\varphi (M^T   (A\widehat{x} -y)) \leq  \varepsilon. $  For this pair $(y, \varepsilon),$ we consider the
problem (\ref{GP}) to which, by the assumption,  any   feasible point   can be approximated by an optimal solution of (\ref{GP}) with   error (\ref{GERROR}).  Therefore, there is an
optimal solution $x^*$ of  (\ref{GP})  which approximates $\widehat{x}$ with   error
 $ \|\widehat{x}- x^*\|_2 \leq C' \sigma_k(\widehat{x})_1 + C'' \gamma (\varepsilon).$
Since $\widehat{x} $ is $k$-sparse, we have $\sigma_k(\widehat{x})_1 =0. $
  Thus  $ \|\widehat{x}- x^*\| _2 \leq  C'' \gamma (\varepsilon) $  which, together with  (\ref{choice}),   implies that for positive components $\widehat{x}_i>0, $ the corresponding components $x^*_i$ must be  positive,  and that for negative components $\widehat{x}_i<0$, the corresponding $x^*_i$ must be  negative. Thus, we have
   \begin{equation}
\label{SET-B} S_1 =\{i: \widehat{x}_i>0 \} \subseteq \{i: x^*_i >0\}   , ~ S_2 =\{i: \widehat{x}_i<0\} \subseteq \{i: x^*_i <0\} . \end{equation}  Note
that $x^*$ is an optimal solution to the convex optimization problem
(\ref{GP}), which satisfies a constraint qualification.  In fact, if $\varepsilon =0,$  the constraint is reduced to the linear equation $Ax=y,$ and if $\varepsilon >0,$ the Slater's constraint qualification is satisfied due to the fact that $A$ is underdetermined (in which case there is a vector $z$ such that $Az=y, $ and hence $\varphi(M^T(Az-y)) <\varepsilon $).   So $x^*$  must satisfy the KKT optimality condition, i.e.,
$$ 0\in \partial_x \left\{\|x \|_1 + \lambda \left[\varphi(M^T(Ax -y)) -\varepsilon\right]     \right\}{\Big |}_{x=x^*},   $$ where $\lambda $ is a Lagrangian multiplier and $\partial_x$ denotes the subgradient  with respect to $x.$  Note that the domains of the functions  $ \varphi(M^T(Ax-y)) $ and $\|x\|_1 $ are $\mathbb{R}^n$, by Theorem 23.8 in \cite{R70},  the above optimality condition is equivalent to  $$   0 \in \left\{ \partial \|x^*\|_1+ \lambda A^TM \partial \varphi(M^T(Ax^*-y))\right\},    $$ where $\partial \varphi(M^T(Ax^*-y))$ is the subgradient of $\varphi$ at $M^T(Ax^*-y)$, and
$ \partial \|x^*\|_1 $ is the subgradient of the $\ell_1$-norm at
$x^*,$ i.e.,    $$ \partial \|x^*\|_1 =\{\zeta \in \mathbb{R}^n:  ~ \zeta_i
=1 \textrm{ for } x^*_i > 0, ~ \zeta_i =-1 \textrm{ for } x^*_i <0, ~
|\zeta_i| \leq   1 \textrm{ for }  x^*_i =0\}. $$
Thus there exists a vector $ v\in \partial \varphi(M^T(Ax^*-y)) $ and a vector  $\zeta \in  \partial \|x^*\|_1 $ such that $ \zeta + \lambda A^T M v    =0.$
 Setting $u= - \lambda Mv $ yields  $\zeta = A^T u.$  Since $\zeta \in  \partial \|x^*\|_1,$ we see that
$    \zeta_i=1$ for $x^*_i>0$, $ \zeta_i=-1$
for $x^*_i< 0,$ and $ |\zeta_i | \leq 1$ for   $x^*_i=0.$  This,
together with (\ref{SET-B}), implies that the vector $\zeta $
satisfies (\ref{RSP}). Note that $S_1$ and $S_2$ are arbitrary and disjoint subsets of $\{1,..., n\}$ with $|S_1|+|S_2| \leq k. $   By Definition 2.6,  $A^T$ must satisfy the weak RSP of order $k.$  ~ $\Box.$

\vskip 0.07in

 As shown by Theorem 2.7,  the weak RSP of $A^T$  is a  necessary condition for many sparse recovery problems to be  stable. It implies that this condition cannot be further relaxed in order to  ensure the stability of these problems. In later sections, we show  that the weak RSP of $A^T$  is also a sufficient condition for a wide range of sparse optimization problems, including  DS and LASSO, to be   stable in sparse data recovery.
 Under the assumptions of RIP, stable NSP or Robust NSP, the traditional error bounds of a recovery problem are measured in terms of these conventional matrix constants. Different from these assumptions, the weak RSP of $A^T$ is a constant-free condition in the sense that the definition of this property does not involve any constant, so is the standard NSP of order $k$ and RSP of order $k.$ Thus an immediately question arises: How to establish the stability of DS and LASSO under a constant-free matrix property?
 The main purpose of this study is to address this question.

\section{Dantzig selectors with linear constraints}

The constraint of (\ref{L1}) becomes linear when $\phi$ is the $\ell_\infty$-norm, $\ell_1$-norm, or their combination. In this case,  problem (\ref{L1}) is equivalent to a linear program  which can be  solved efficiently by simplex methods or interior-point methods. Thus in this section, we focus on the norm   $\phi(\cdot) = \alpha \| \cdot\|_\infty + (1-\alpha) \| \cdot \|_1, $  where $\alpha \in [0, 1]$ is a  fixed constant. Using this norm,  the DS is in the form
\begin{equation}  \label{GDS}
\min_{x} \{\|x\|_1: ~ \alpha \|M^T(Ax-y)\|_\infty +(1-\alpha) \|M^T(Ax-y)\|_1 \leq \tau\}.
\end{equation}
 This problem  encompasses the following special cases:
\begin{equation}  \label{DS}
\min_{x} \{\|x\|_1:  ~  \|M^T(Ax-y)\|_\infty  \leq \tau\},
\end{equation}
\begin{equation}  \label{DS1}
\min_{x} \{\|x\|_1:   ~ \|M^T(Ax-y)\|_1 \leq \tau\},
\end{equation}
which correspond  to the cases  $\alpha =1$ and $\alpha =0 $ in (\ref{GDS}), respectively.  Particularly, when $M=A,$  problem (\ref{DS})  is   the   standard DS  proposed by Cand\`es and Tao \cite{CT07}.
The purpose of this section is to  establish a  stability result for (\ref{GDS}) under the weak RSP of order $k$ of $A^T.$

By introducing two auxiliary variables $(\xi, v)$, the DS problem (\ref{GDS})  can be written as
\begin{eqnarray}    \label{GDS-01}   & \min_{(x, \xi, v)}   & \|x\|_1   \nonumber\\
& \textrm{s.t.} &   \alpha \xi +(1-\alpha) \mathbf{e}^T v \leq \tau,  ~  \|M^T(Ax-y)\|_\infty  \leq \xi,  \\
 & &   |M^T(Ax-y)| \leq v,  ~ \xi \in  \mathbb{R}_+,  ~ v \in  \mathbb{R}^q_+ ,  \nonumber
 \end{eqnarray} where $ \mathbf{e}$ is the vector of ones in $\mathbb{R}^q. $
  Introducing the variable  $t \in  \mathbb{R}^n_+$ such that $|x|\leq t$,   the   problem (\ref{GDS-01})
 can be   written as
\begin{eqnarray}  \label{GDS-02}
& \min_{(x, t, \xi, v)}   & \widetilde{\mathbf{e}}^T t  \nonumber \\
& \textrm{s.t.} &   x\leq t,  ~-x\leq t,  ~ \alpha \xi +(1-\alpha) \mathbf{e}^T v \leq \tau,  \nonumber \\
 & &    -\xi \mathbf{e}   \leq  M^T(Ax-y) \leq  \xi \mathbf{e}, ~ M^T(Ax-y) \leq v, \\
& &   -M^T(Ax-y) \leq v,  ~ t \in  \mathbb{R}^n_+, ~ \xi \in  \mathbb{R}_+,  ~ v \in  \mathbb{R}^q_+,  \nonumber
\end{eqnarray} where $\widetilde{\mathbf{e}} $ is the vector of ones in $ \mathbb{R}^n. $
Clearly, at any optimal solution $(x, t, \xi, v) $ of   (\ref{GDS-02}), it must hold that $t=|x|, ~\xi\geq  \|M^T(Ax-y)\|_\infty $ and $ v \geq |M^T (Ax-y)|. $
There are many ways to obtain the dual problem of (\ref{GDS-02}).  For example, we may rewrite (\ref{GDS-02}) as the so-called canonical form  $\min_{z} \{r^T z:  ~Qz\geq \nu, ~z \geq 0\} $ with problem data $(Q, \nu, r) ,$ to which
the Lagrangian dual is given as   $\max_w \{\nu^Tw:  ~Q^T w\leq r, ~w\geq 0\}. $  So it is   easy to verify that the  dual problem of  (\ref{GDS-02}) is given as
\begin{eqnarray}  \label{DUAL-DS}
& \max_{(w^{(i)}, i=1, \dots , 7)}   &  -\tau w^{(3)} +y^T M (w^{(4)}-w^{(5)}- w^{(6)} +w^{(7)})   \nonumber \\
& \textrm{s.t.} & w^{(1)}+w^{(2)} \leq  \widetilde{\mathbf{e}},
   \nonumber\\
    & &  -w^{(1)} + w^{(2)}+A^TM  (w^{(4)}-w^{(5)}-w^{(6)}+ w^{(7)}) =0,  \\
     & & -\alpha w^{(3)} +\mathbf{e}^Tw^{(4)}+ \mathbf{e}^Tw^{(5)}  \leq 0,  \nonumber \\
     & & -(1-\alpha) w^{(3)} \mathbf{e}  +  w^{(6)} +w^{(7)}   \leq 0,   \nonumber\\
    & & w^{(i)}\geq 0, ~ i=1, \dots , 7, \nonumber
\end{eqnarray}
where $w^{(1)}, w^{(2)}\in \mathbb{R}^n_+, w^3\in \mathbb{R}_+$ and $ w^{(4)}, \dots, w^{(7)} \in \mathbb{R}^q_+. $
By the KKT optimality condition of (\ref{GDS-02}) or (\ref{DUAL-DS}), we immediately have the following lemma.

\vskip 0.07in

  \textbf{Lemma 3.1.}    \emph{$x^*$ is an optimal solution of the DS problem
(\ref{GDS}) if and only if there exist vectors $ (t^*,  \xi^*,  v^*, w^{(1)}_*, \dots, w^{(7)}_*) $
such that $ (x^*, t^*, \xi^*, v^*, w^{(1)}_*, \dots, w^{(7)}_* ) \in \mathfrak{D},  $
where $\mathfrak{D}$ is the set of vectors  $(x,t,\xi,v,w^{(1)},\dots, w^{(7)})$ satisfying that
\begin{equation}\label{KKT} \left\{\begin{array} {l}
    x\leq t,  ~ -x\leq t, ~ \alpha \xi +(1-\alpha) \mathbf{e}^T v \leq \tau,  \\
    -\xi \mathbf{e}   \leq  M^T(Ax-y) \leq  \xi \mathbf{e},  ~ M^T(Ax-y) \leq v, ~ -M^T(Ax-y) \leq v,   \\
  -w^{(1)} + w^{(2)}+A^TM  (w^{(4)}-w^{(5)}-w^{(6)} + w^{(7)}) =0,
  \\
      w^{(1)}+w^{(2)} \leq \widetilde{\mathbf{e}}, ~ -\alpha w^{(3)} +\mathbf{e}^Tw^{(4)}+ \mathbf{e}^Tw^{(5)}  \leq 0,     \\
       -(1-\alpha) w^{(3)} \mathbf{e}+  w^{(6)} +  w^{(7)}   \leq 0,  \\
      -\tau w^{(3)} +y^T M (w^{(4)}-w^{(5)}-w^{(6)}+ w^{(7)}) =\widetilde{\mathbf{e}}^T t, \\
       t \in  \mathbb{R}^n_+,  ~  \xi \in  \mathbb{R}_+,  ~ v \in  \mathbb{R}^q_+ , ~  w^{(i)}\geq 0, ~ i=1, \dots , 7.
    \end{array} \right. \end{equation}
 For any $(x,t,\xi,v, w^{(1)}, \dots, w^{(7)} ) \in \mathfrak{D},  $ we must have that $t=|x|, ~ \xi\geq  \|M^T(Ax-y)\|_\infty $ and $ v \geq |M^T (Ax-y)|. $}

Note that the conditions $t \in  \mathbb{R}^n_+,  ~  \xi \in  \mathbb{R}_+ $ and $ v \in  \mathbb{R}^q_+ $ are implied from other conditions in (\ref{KKT}), and thus these conditions can be removed from (\ref{KKT}).  It is easy to write   $\mathfrak{D}$ in the form
\begin{equation} \label{D} \mathfrak{D} = \{z=(x,t,\xi,v, w^{(1)},\dots, w^{(7)} ):  ~ M^1 z  \leq d^1 ,  ~ M^2
z = d^2 \}, \end{equation}  where $d^2=0, $  $d^1= \left(0, 0, \tau,y^TM , - y^TM,y^TM, -y^TM,\widetilde{\mathbf{e}}^T ,0 ,0 ,0 ,0 , 0 ,0 ,0 ,0 ,0
                      \right)^T$  and
{\small \begin{equation} \label{MM-X}
 M^1=  \left[
        \begin{array}{ccccccccccc}
         I  & -I & 0 & 0 & 0 & 0 & 0 & 0 & 0 & 0 & 0 \\
        -I  & -I & 0 & 0 & 0 & 0 & 0 & 0 & 0 & 0 & 0 \\
          0 & 0 & \alpha & (1-\alpha)\mathbf{e}^T & 0 & 0 & 0 & 0 & 0 & 0 & 0 \\
         M^TA & 0 & -\mathbf{e} & 0 & 0 & 0 & 0& 0 & 0 & 0 & 0 \\
          -M^TA & 0 & -\mathbf{e} & 0 & 0 & 0 & 0 & 0 & 0 & 0 & 0 \\
          M^TA & 0 & 0 & -I & 0 & 0 & 0 & 0 & 0 & 0 & 0 \\
          -M^TA & 0 & 0 & -I & 0 & 0 & 0 &0 & 0 & 0 & 0 \\
          0 & 0 & 0 & 0 & I & I & 0 & 0 & 0 & 0 & 0 \\
          0 & 0 & 0 & 0 & 0 & 0 & -\alpha & \mathbf{e}^T & \mathbf{e}^T & 0 & 0 \\
          0 & 0 & 0 & 0 & 0 & 0 & -(1-\alpha)\mathbf{e} & 0& 0 & I & I \\
          0 & 0 & 0 & 0 & -I & 0 & 0 & 0 & 0 & 0 & 0 \\
          0 & 0 & 0 & 0 & 0 & -I & 0 & 0 & 0 & 0 & 0 \\
          0 & 0 & 0 & 0 & 0 & 0 & -1 & 0 & 0 & 0 & 0 \\
          0 & 0 & 0 & 0 & 0 & 0 & 0 & -I & 0 & 0 & 0 \\
          0 & 0 & 0 & 0 & 0 & 0 & 0 & 0 & -I & 0 & 0 \\
          0 & 0 & 0 & 0 & 0 & 0 & 0 & 0 & 0 & -I & 0 \\
          0 & 0 & 0 & 0 & 0 & 0 & 0 & 0 & 0 & 0 & -I \\
        \end{array}
      \right], \end{equation} }
 \begin{equation} \label{MM-Y}
         M^2 = \left[
           \begin{array}{ccccccccccc}
             0 & 0 & 0 & 0 & -I & I & 0 &  A^TM  & -A^TM & -A^TM   & A^TM\\
             0 & \widetilde{\mathbf{e}}^T & 0 & 0 & 0 & 0 & \tau &  -y^TM  & y^TM & y^TM   & -y^TM
               \end{array}
               \right], \end{equation}
               where $I $'s are the  identity matrices and 0's are the zero matrices with   suitable sizes.
We now prove that the stability  of (\ref{GDS}) is guaranteed when $A^T$ admits the weak RSP of order
$k. $

 \vskip 0.07in

\textbf{Theorem 3.2.}    \emph{Let the problem  data $(A, M, y,
\tau, \alpha)$  in (\ref{GDS}) be given, where   $A\in
\mathbb{R}^{m\times n}~  (m<n) $ and $M\in
\mathbb{R}^{m\times q}~  (m\leq q) $ with $\textrm{rank} (A) =\textrm{rank} (M) =m.$
Suppose that $A^T$ satisfies the weak RSP of order $k. $
   Then for any $x \in \mathbb{R}^n, $
    there is an optimal solution $x^*$ of (\ref{GDS})   approximating $x$ with error
   \begin{equation} \label{INQ-AA}
    \| x   -
   x^*   \|_2
\leq    \gamma   \left\{\left(\alpha \|\vartheta\|_\infty +(1-\alpha) \|\vartheta\|_1 - \tau\right)^+ + 2\sigma_k (x)_1  +  c(\tau     +    \|\vartheta\|_\infty  )      \right\},
\end{equation}
 where $\vartheta =M^T (Ax-y),$  $c  $ is the constant given as   $$c= \max_{G\subseteq \{1,\dots,q\}, |G | =m} \{ \|M_G^{-1} (AA^T)^{-1} A \|_{\infty\to 1}:   M_G \in \mathbb{R}^{m\times m} \textrm{  is an invertible submatrix of }  M \}, $$   and $\gamma
=\sigma_{\infty,2} (M^1,M^2) $
is the Robinson's constant with  $(M^1,M^2)$ being given
in (\ref{MM-X}) and (\ref{MM-Y}).  In particular, for any $x$ satisfying the constraints of (\ref{GDS})
     there is a solution $x^*$ of (\ref{GDS})   approximating  $x$ with error
         \begin{equation} \label{INQ-new}
    \| x   -
   x^* \|_2
\leq    \gamma   \left\{  2\sigma_k (x)_1  +  c\tau     +    c\| M^T(Ax-y) \|_\infty      \right\}  \leq 2\gamma   \left\{  \sigma_k (x)_1  +  c \tau     \right\}.
\end{equation} }

\emph{Proof.}  Let $x$ be any given vector in $ \mathbb{R}^n. $ We set  $(t, \xi, v)$   as follows:
\begin{equation} \label{ttvv} t=|x|,  ~  \xi =\|M^T(Ax-y) \|_\infty, ~  v= |M^T(Ax-y)|.   \end{equation}
   Let $S$ be the support of the $k$-largest components of
 $|x|,$ i.e., $S= \{i_1, \dots, i_k\},$ when the components of $|x|$ are sorted into a descending order as $|x_{i_1}|\geq \cdots \geq |x_{i_k}| \geq  \cdots \geq |x_{i_n}|. $ We define by $S' =\{j\in S:  x_j> 0\} $
 and $ S'' = \{j\in S: x_j <  0\}  $  which imply that $S'$ and $ S''$ are disjoint and with cardinality $|S' \cup S''|=|S| \leq k. $
Note that $A^T$ satisfies the weak RSP of order $k.$    There exists a vector   $  \zeta=  A^T u^*   $
 for some $ u^*  \in \mathbb{R}^m$  satisfying the following properties:
  $  \zeta_i= 1\textrm{ for }i\in S', ~ \zeta_i=-1 \textrm{ for   }i\in S'', $ and
 $ |\zeta_i| \leq 1\textrm{ for  } i \notin S' \cup S'' .$

We now construct a
vector $(w^{(1)}, \dots, w ^{(7)})$
which satisfies the constraint of  problem (\ref{DUAL-DS}).
First, we set $ w^{(1)} $ and $ w^{(2)}$ as follows:  $ w^{(1)}_i= 1\textrm{ and } w^{(2)}_i=0\textrm{ for
}i\in S'; $    $ w^{(1)}_i=0\textrm{ and  } w^{(2)}_i=
1 \textrm{ for }i\in S'';  $   and $ w^{(1)}_i =(1+\zeta_i)/2 $ and
$w^{(2)}_i = (1-\zeta_i)/2 \textrm{ for all }i\notin S'\cup S'' .
$  This choice of $(w^{(1)}, w^{(2)})$ satisfies that $w^{(1)} - w^{(2)}
=\zeta.$  Note that $M$ is a full-row-rank matrix. Thus there exists an invertible $m\times m$ submatrix of $M$, denoted by $M_\mathfrak{J}$, where $\mathfrak{J}\subseteq \{1, \dots, q\}$ with cardinality $|\mathfrak{J}|=m.$  We choose the vector $h \in \mathbb{R}^q $ as
$h_{\mathfrak{J}} = M^{-1}_{\mathfrak{J}}  u^*  $ and $h_{\overline{\mathfrak{J}}} =0  $ where  $ \overline{\mathfrak{J}} =\{1, \dots , q\}\backslash \mathfrak{J}. $  Thus $M h=u^*.$
We  now define the nonnegative vectors $w^{(3)}, \dots, w^{(7)} $ as follows:
$$ w^{(3)} = \|h\|_1,  ~ w^{(4)} = \alpha (h)^+,  ~ w^{(5)} = -\alpha (h)^-, ~ w^{(6)} =-(1-\alpha) (h)^- , ~ w^{(7)} = (1-\alpha) (h)^+, $$
which implies that
$$ \mathbf{e}^T(w^{(4)}+w^{(5)})  = \alpha \mathbf{e}^T |h|= \alpha \|h\|_1= \alpha w^{(3)},
 ~  w^{(6)} + w^{(7)}   = (1-\alpha) |h|     \leq  (1-\alpha) w^{(3)} \mathbf{e}.$$
\begin{equation} \label{W4567}  w^{(4)}-w^{(5)}- w^{(6)} + w^{(7)} = \alpha (h)^+ + \alpha (h)^- +(1-\alpha) (h)^- + (1-\alpha) (h)^+   =h . \end{equation}
 Therefore,
$$ A^T M  (w^{(4)}-w^{(5)}- w^{(6)} + w^{(7)} ) = A^T (M h) =  A^Tu^*  = \zeta = w^{(1)}-w^{(2)}.$$  Thus the  vector $(w^{(1)}, \dots, w ^{(7)})$ constructed as above satisfies the constraints of   (\ref{DUAL-DS}).
Consider the set  $\mathfrak{D}$   written as (\ref{D}). For the vector $(x, t, \xi,   v, w^{(1)}, \dots , w ^{(7)} ), $ where $(t, \xi,v)$ is chosen as (\ref{ttvv}) and $ w^{(1)}, \dots , w ^{(7)}$ are chosen as above, applying Lemma 2.1 with $(M^1,M^2)
   $ being given in (\ref{MM-X}) and (\ref{MM-Y}),
   there  exists  a vector
  $(x^*,t^*,\xi^*, v^*,  w^{(1)}_*, \dots , w ^{(7)}_*)  \in
   \mathfrak{D}  $ such that
\begin{equation} \label {hoff} \left\|
   \left[
      \begin{array}{c}
        x \\
        t \\
       \xi  \\
        v \\
        w^{(1)} \\
        \vdots \\
        w^{(7)} \\
      \end{array}\right] - \left[
       \begin{array}{c}
        x^* \\
         t^* \\
        \xi^* \\
         v^* \\
        w^{(1)}_* \\
        \vdots \\
         w^{(7)}_* \\
       \end{array}
     \right]
 \right\|_2 \leq  \gamma  \left\| \left[\begin{array}{c}
   (x-t)^+\\
   (-x- t)^+ \\
  (\alpha \xi +(1-\alpha) \mathbf{e}^T v - \tau)^+ \\
   (M^T(Ax-y)-\xi \mathbf{e})^+ \\
   (-M^T(Ax-y)-\xi \mathbf{e})^+\\
   (M^T(Ax-y)-v)^+ \\
   (-M^T(Ax-y)-v)^+\\
  A^TM  (w^{(4)}-w^{(5)}-w^{(6)} + w^{(7)}) -w^{(1)} + w^{(2)} \\
      (w^{(1)}+w^{(2)} - \widetilde{\mathbf{e}})^+   \\
       ( -\alpha w^{(3)} +\mathbf{e}^Tw^{(4)}+ \mathbf{e}^Tw^{(5)} )^+      \\
       (-(1-\alpha) w^{(3)} \mathbf{e}+  w^{(6)} +  w^{(7)}  )^+    \\
      \widetilde{\mathbf{e}}^T t + \tau w^{(3)} - y^T M (w^{(4)} - w^{(5)} - w^{(6)}+ w^{(7)})   \\
      (Y)^+
            \end{array}
             \right] \right\|_1, \end{equation}
where $ (Y)^+ = ( (-t)^+, (-\xi)^+,  (-v)^+,
  (-w^{(1)})^+, \dots, (-w^{(7)})^+), $  and $\gamma  =$$\sigma_{\infty, 2}
(M^1, M^2) $  is the Robinson's constant with $ (M^1,
M^2) $ being given in  (\ref{MM-X}) and (\ref{MM-Y}).    It follows from (\ref{ttvv})  that
$$  (x-t)^+=(-x- t)^+ =0 ,  \left(M^T(Ax-y)-\xi \mathbf{e}\right)^+ =0 , $$  $$\left(-M^T(Ax-y)-\xi \mathbf{e}\right)^+ =0,  ~ \left(M^T(Ax-y)-v\right)^+  =0, ~\left(-M^T(Ax-y)-v\right)^+ =0 .  $$
 Since $ (w^{(1)}  ,\dots ,  w^{(7)}  )$  is a  feasible point to (\ref{DUAL-DS}), we also see that
$$  A^TM  (w^{(4)}-w^{(5)}-w^{(6)} + w^{(7)})-w^{(1)} + w^{(2)}=0,  ~ \left(w^{(1)}+w^{(2)} - \widetilde{\mathbf{e}}\right)^+ =0,$$
    $$   \left( -\alpha w^{(3)} +\mathbf{e}^Tw^{(4)}+ \mathbf{e}^Tw^{(5)} \right)^+   =0,
     ~   \left(-(1-\alpha) w^{(3)} \mathbf{e}+  w^{(6)} +  w^{(7)}  \right)^+  =0 . $$  Moreover, the nonnegativity of  $(t, \xi, v, w^{(1)}, \dots , w ^{(7)})$  implies that  $ (Y)^+  =0 .$     Note that  $$   \|  x   -
   x^*  \|_2 \leq \|(x ,t ,\xi , v , w^{(1)} , \dots, w^{(7)}) - (
        x^*,
         t^*,
        \xi^*,
         v^*,
        w^{(1)}_*, \dots,
         w^{(7)}_* )\|_2. $$  Thus the inequality (\ref{hoff}) is reduced to
\begin{equation} \label{INQ-22}
      \|  x   -  x^*  \|_2  \leq    \gamma  \left\| \left[\begin{array}{c}
  \left(\alpha \xi +(1-\alpha) \mathbf{e}^T v - \tau\right)^+ \\
      \widetilde{\mathbf{e}}^T t +\tau w^{(3)} - y^T M (w^{(4)}-w^{(5)}-w^{(6)}+ w^{(7)})
            \end{array}
             \right] \right\|_1.
\end{equation}
 Denote by $\vartheta =M^T(Ax-y) $ which implies that $y^T M = x^T A^T M -\vartheta ^T.$
   It follows from (\ref{ttvv}), (\ref{W4567}) and the fact $A^TMh=A^Tu^*=\zeta$  that
   \begin{equation} \label{F1} \left(\alpha \xi +(1-\alpha) \mathbf{e}^T v - \tau\right)^+ = \left(\alpha \|\vartheta\|_\infty +(1-\alpha) \|\vartheta\|_1 - \tau\right)^+  \end{equation}
   and
\begin{eqnarray} \label{F2}     & &    \left |   \widetilde{\mathbf{e}}^T t + \tau w^{(3)} - y^T M (w^{(4)}-w^{(5)}-w^{(6)}+ w^{(7)})  \right |
 \nonumber \\ &     &  =   \left |  \widetilde{\mathbf{e}}^T|x|  + \tau w^{(3)} - (x^T A^T M-\vartheta^T)  h  \right |
     =     \left |  \widetilde{\mathbf{e}}^T|x|  + \tau w^{(3)} -  x^T \zeta + \vartheta^T   h  \right |   \nonumber  \\
 &   & \leq         2\sigma_k (x)_1  + \tau \|h\|_1   +   \|\vartheta\|_\infty   \|h\|_1,
\end{eqnarray}
where the last inequality follows from the fact $ |\vartheta^Th| \leq  \|\vartheta\|_\infty   \|h\|_1 $ and  $$ \left|  \widetilde{\mathbf{e}}^T|x|    -  x^T \zeta \right|= \left|\widetilde{\mathbf{e}}^T|x|-  x^T_S \zeta_S - x^T_{\overline{S}} \zeta_{\overline{S}}\right|  \leq  \|x\|_1- \| x _S \|_1 +  \|x^T_{\overline{S}}\|_1 \|\zeta_{\overline{S}}\|_\infty \leq 2 \|x^T_{\overline{S}}\|_1 =2\sigma_k(x)_1. $$
We define the constant
$$ c= \max_{G\subseteq \{1,\dots, q\}, |G|=m} \{ \|M_G^{-1} (AA^T)^{-1} A \|_{\infty\to 1}:   M_G \in \mathbb{R}^{m\times m} \textrm{  is an invertible submatrix of }  M \}. $$  Noting that  $h_{\overline{\mathfrak{J}}}=0 $ and $M_\mathfrak{J} h_\mathfrak{J} = u^*= (AA^T)^{-1}A \zeta,$ we have
\begin{equation}  \label{F3}  \|h\|_1=   \|h_\mathfrak{J}\|_1 = \|M_\mathfrak{J}^{-1} (AA^T)^{-1} A \zeta \|_1 \leq \|M_\mathfrak{J}^{-1} (AA^T)^{-1} A\|_{\infty \to 1} \|\zeta\|_\infty  \leq c. \end{equation}
Combining  (\ref{INQ-22})--(\ref{F3})   leads to
$$
    \| x   -
   x^*   \|_2
\leq    \gamma   \left\{\left(\alpha \|\vartheta\|_\infty +(1-\alpha) \|\vartheta\|_1 - \tau\right)^+ + 2\sigma_k (x)_1  +  c(\tau     +    \|\vartheta\|_\infty  )        \right\},
$$
which is exactly the bound given in (\ref{INQ-AA}). In particular, if $x$ satisfies the constraint of (\ref{GDS}), then
\begin{eqnarray} \label{INQ-2}
    \| x   -
   x^*   \|_2
\leq    \gamma   \left\{  2\sigma_k (x)_1  +  c (\tau     +    \|\vartheta\|_\infty  )       \right\}.
\end{eqnarray}
Since $ \|\vartheta\|_\infty \leq \alpha \|\vartheta\|_\infty + (1-\alpha) \|\vartheta\|_1 \leq \tau,$  the  bound (\ref{INQ-new}) follows from (\ref{INQ-2}) immediately.  ~~ $ \Box $

\vskip 0.07in

Many existing conditions imply the weak RSP of order $k$ of $A^T.$  The following
 result can be immediately obtained from Theorem 3.2.

\vskip 0.07in

\textbf{Corollary 3.3.} \emph{Let $A$ and $M$ be given as in Theorem 3.2. Suppose that one of the following conditions holds:}
  (a)  \emph{$A$ (with $\ell_2$-normalized columns) satisfies  the mutual coherence property
$\mu_1(k) + \mu_1(k-1) <1; $}   (b)  \emph{RIP of order $2k$ with constant $\delta_{2k} < 1/\sqrt{2};  $  } (c)  \emph{stable NSP of order $k$ with constant $\rho \in(0,1); $}
(d)  \emph{robust NSP of order $k$  with  constants $ \rho' \in (0,1)$
and $\rho'' >0; $}
(e)   \emph{NSP of order $k; $ }
(f)  \emph{RSP of order $k $ of $A^T.$ Then the conclusions of Theorem 3.2 are valid for  (\ref{GDS}).}

 \vskip 0.07in

In this corollary,  condition (a)  implies (e) (see Theorem 5.15 in \cite{EK12} and the definition of mutual coherence $\mu_1(\cdot)$ therein). Condition (b) implies that every $k$-sparse vector can be exactly recovered by $\ell_1$-minimization (see \cite{CZ14}), and thus implies (e). Each of conditions (c) and (d) implies  (e). Conditions (e) and (f) are equivalent.  Thus each of conditions (a) -- (e) implies (f), and hence they imply the weak RSP of order $k$ of $A^T.$  Thus Corollary 3.3 follows  immediately from Theorem 3.2.
This is a unified result   in the sense that any of the  above-mentioned conditions  implies the same error bounds (\ref{INQ-AA}) and (\ref{INQ-new}).
     In particular, by setting $\alpha =1 $ and $ \alpha =0,$  respectively,   Theorem 3.2 claims that under the weak RSP of order $k$, both   DS problems (\ref{DS}) and (\ref{DS1}) are   stable in  sparse vector recovery. The results  for   $M=I$ and $M=A,$ respectively, can  follow immediately from Theorem 3.2 and Corollary 3.3.   As an example, let us state the result for the standard DS, corresponding to the case $M=A$ and $\alpha=1$ in (\ref{GDS}). Combining Theorems 3.2 and 2.7 as well as  Corollary 3.3 leads to the following result.

\vskip 0.07in

  \textbf{Corollary 3.4.}   \emph{Let  $A\in
\mathbb{R}^{m\times n} $ with $\textrm{rank} (A) =m < n.$  Consider the following standard  DS problem: \begin{equation}  \label{CT-DS}  \min_x\{\|x\|_1: \|A^T(Ax-y)\|_\infty  \leq \tau\}. \end{equation}   Then the  following statements  hold:
\begin{itemize}
\item[(i)] Suppose that $A^T$ satisfies the weak RSP of order $k. $
  Then for any $x \in \mathbb{R}^n, $
    there is an optimal solution $x^*$ of   (\ref{CT-DS}) approximating $x$ with error
 $$
    \| x   -
   x^*   \|_2
\leq    \gamma   \left\{(  \| A^T (Ax-y)\|_\infty   - \tau)^+ + 2\sigma_k (x)_1  +  c_A \tau     +    c_A \|A^T (Ax-y)\|_\infty        \right\},
$$
 where   $c_A $ is the constant  $$c_A= \max_{G\subseteq \{1,\dots,n\}, |G | =m} \left\{ \|A_G^{-1} (AA^T)^{-1} A \|_{\infty\to 1}:   A_G   \textrm{ is an invertible submatrix of }  A \right\}, $$   and $\gamma
=\sigma_{\infty,2} (M^1,M^2) $
is the Robinson's constant determined by  $(M^1,M^2)$ given
in (\ref{MM-X}) and (\ref{MM-Y}) with $M=A$ and $\alpha=1.$  In particular, for any $x$ satisfying
the constraint of (\ref{CT-DS}),
     there is a solution $x^*$ of   (\ref{CT-DS})  approximating $x$ with error
         \begin{equation} \label{E2}
    \| x   -
   x^* \|_2
\leq    \gamma   \left\{  2\sigma_k (x)_1  +  c_A \tau     +   c_A \| A^T (Ax-y) \|_\infty       \right\}  \leq 2\gamma   \left\{  \sigma_k (x)_1  +  c_A \tau     \right\}.
\end{equation}
  Conversely, if for any given small $\tau \geq 0 $ and for any $x$ obeying $  \|A^T(Ax-y)\|_\infty    \leq \tau, $ there is an optimal solution $x^*$ of (\ref{CT-DS}) such that the estimate (\ref{E2}) holds,  then $A^T$  must satisfy the weak RSP of order $k.  $
\item [(ii)]  Every matrix condition  listed in Corollary 3.3 is sufficient to ensure that, for any $x$ obeying $  \|A^T(Ax-y)\|_\infty    \leq \tau, $ there is an optimal solution $x^*$ of (\ref{CT-DS}) such that  (\ref{E2}) holds.
 \end{itemize} }

 \vskip 0.07in

    Item (i)   above is   different from existing stability results for the standard DS in terms of the mild assumption, analytic method, and the way of expression of stability coefficients. Roughly speaking, Corollary 3.4 indicates that the weak RSP of $A^T$ is  both necessary and sufficient for the standard DS  to be weakly stable in sparse recovery. Item (ii) above indicates that  the error bound   (\ref{E2})   holds under any of the conditions listed in Corollary 3.3. Letting $\alpha=0$ and replacing $\|\cdot\|_\infty $ by $\|\cdot\|_1, $ Corollary 3.4 immediately becomes the stability result for the problem (\ref{DS1}) with $M=A.$




  \vskip 0.07in

  \textbf{Remark 3.5.}   When the matrix $A$ does not satisfy a desired matrix property like the RIP, NSP, RSP or REC), a scaled version of this matrix, i.e., $AU,$ where $U$ is a nonsingular matrix, might admit a desired property.  This partially explains why a weighted $\ell_1$-minimization algorithm  (e.g., \cite{CWB08, ZL12, ZK15})  often numerically outperforms  the standard  $\ell_1$-minimization in sparse data recovery.    The  stability theory developed in this paper can be easily generalized to  weighted $\ell_1$-minimization  problems. Take the following weighted Dantzig selector  as an example:
  \begin{equation} \label{WDS}  \min \{\|Wx\|_1:  ~ \|A^T(Ax-y)\|_\infty \leq \tau\}, \end{equation}   where $W$ is a nonsingular diagonal matrix.
We ask whether this problem is   weakly stable in sparse data recovery. By the nonsingular transform $u=Wx,$ this problem can be written as
   $$ \min \{\|u\|_1:  ~ \|A^T(AW^{-1} u-y)\|_\infty \leq \tau\},$$
which is of the form
\begin{equation}\label{ww-DS}  \min \{\|u\|_1:  ~ \|M^T(\widetilde{A} u-y)\|_\infty \leq \tau\},\end{equation}
with $ \widetilde{A} =AW^{-1}$ and $M=A. $ This is   the recovery problem (\ref{L1}) with $\phi=\|\cdot\|_\infty. $   Thus it is straightforward to extend the stability results developed in this paper to the  weighted Dantzig selector (\ref{WDS}) under the weak RSP assumption on the scaled matrix  $ \widetilde{A} =A W^{-1}. $

\section{Dantzig selectors with nonlinear constraints}

In this section,  we deal with the nonlinear problem (\ref{L1}), where
  the constraint $\phi(M^T(Ax-y)) \leq \tau  $  cannot be represented exactly as a finite number of linear constraints, for example, when $\phi =\|\cdot \|_p$ with $p\in (1, \infty).$  In this case, $\tau$ must be positive, since otherwise if $\tau=0$ the constraint will reduce to the linear system $Ax=y.$  We show that the nonlinear DS problem (\ref{L1}) remains stable in sparse data recovery under the weak RSP assumption.

   Let $\varrho^* $ be  the optimal value  of  (\ref{L1}) and $S^*$ the set of optimal
solutions of (\ref{L1}),
which clearly can be written as
$$ S^* =\{ x\in \mathbb{R}^n : ~ \|x\|_1\leq \varrho^* ,  ~ \phi(M^T(Ax-y)) \leq \tau\}. $$   In terms of   $\mathfrak{B}^\phi $ in $\mathbb{R}^q,$   the nonlinear problem (\ref{L1}) can be written  as
\begin{equation} \label{l1bb} \varrho^* = \min_{(x, u)} \{ \|x\|_1: ~ u=
M^T(Ax-y)/\tau,   ~ u\in \mathfrak{B}^\phi \}.
 \end{equation}
Unlike the linear case  examined  in Section 3, the nonlinearity of the constraint prohibits  applying  Lemma 2.1 directly to establish a stability result. A natural idea is to use a certain polytope approximation of the   unit ball in $\mathbb{R}^q.$     In this section, we use the  polytope $   {\cal Q}_{\epsilon_j}$, which is defined by  (\ref {PJ-01}) and  satisfies the relation (\ref{contain}), to  approximate the unit ball $\mathfrak{B}^\phi .$
Let us first develop a few  technical results. The first one below,  which is of independent interest, describes certain properties of the projection operator.

\vskip 0.07in

\textbf{Lemma 4.1.}  (i) \emph{Let $\Omega \subseteq  T$ be two compact convex
sets in $\mathbb{R}^n. $   Then  for any $x \in \mathbb{R}^n $, \begin{equation}  \label{PPII}
\|\Pi_ \Omega (x) - \Pi_ {T} (x) \|_2^2 \leq d^{\cal H} (\Omega ,
T)  \|x- \Pi_ \Omega (x) \|_2   .
\end{equation}}
(ii)  \emph{ Let $\Omega, U,   T$ be three compact convex
sets in $\mathbb{R}^n  $ satisfying $ \Omega\subseteq  T$ and $ U \subseteq T.$ Then for any $x \in \mathbb{R}^n $ and  any $u\in U$,
\begin{equation}  \label{SUT} \|x- \Pi_ \Omega(x)\|_2 \leq d^{\cal H} (\Omega, T) + 2 \|x-u\|_2. \end{equation}
}
\emph{Proof. } By the property of the projection operator, we have
$
 (x- \Pi_ {T}(x))^T (u-\Pi_ {T}(x)) \leq 0 \textrm{ for all } u \in
T,$ and $(x- \Pi_ \Omega(x))^T (v-\Pi_ \Omega(x)) \leq 0 \textrm{
for all }v \in \Omega. $  Since $ \Pi_ \Omega (x) \in \Omega \subseteq T $ and $\Pi_ \Omega ( \Pi_ T(x))  \in \Omega, $ we immediately have the following two equalities:
 \begin{equation} \label{Proj-bb}  (x- \Pi_ {T}(x))^T [\Pi_ \Omega (x) -\Pi_ {T}(x)] \leq 0,  ~  (x- \Pi_ \Omega(x))^T ( \Pi_ \Omega ( \Pi_ T(x))-\Pi_ \Omega(x)) \leq 0.
\end{equation}
Since  $\Omega\subseteq T $ and $\Pi_ {T}(x) \in T,$ by the definition
  of Hausdorff metric,    we have
 \begin{equation} \label{III} d^{\cal H}
(\Omega , T) = \sup_{w\in T} \inf_{z\in \Omega} \|w-z\|_2 \geq \inf_{z\in \Omega} \|\Pi_ T(x) -z\|_2 = \| \Pi_ {T}(x)- \Pi_ \Omega (\Pi_ T(x))   \|_2.  \end{equation}
 By (\ref{Proj-bb}) and (\ref{III}), we have
  \begin{eqnarray*}
  \|\Pi_ \Omega(x) -\Pi_ {T}(x)\|_2^2  & = & (\Pi_ \Omega(x) -x+x-\Pi_ {T}(x))^T (  \Pi_ \Omega(x) -\Pi_ {T}(x) )\\
   & = &
 (\Pi_ \Omega(x)-x )^T ( \Pi_ \Omega(x)
 -\Pi_ {T}(x)) +  (x- \Pi_ {T}(x))^T ( \Pi_ \Omega(x)
 -\Pi_ {T}(x)) \\
 &\leq  &    (\Pi_ \Omega(x)- x )^T ( \Pi_ \Omega(x)
 -\Pi_ {T}(x))  \\
 & =  &   (\Pi_ \Omega(x)-x )^T [ \Pi_ \Omega(x) - \Pi_ \Omega (\Pi_ T(x)) + \Pi_ \Omega (\Pi_ T(x))
 -\Pi_ {T}(x)]  \\
 &  \leq &   (\Pi_ \Omega(x) - x )^T [\Pi_ \Omega (\Pi_ T(x))
 -\Pi_ {T}(x)]\\
 & \leq & \| x- \Pi_ \Omega(x)\|_2 \|  \Pi_ \Omega (\Pi_ T(x)) - \Pi_ {T}(x) \|_2   \\
 & \leq & d^{\cal H}
(\Omega , T)  \| x- \Pi_ \Omega(x)\|_2 .
 \end{eqnarray*}
Thus (\ref{PPII}) holds. We now prove (\ref{SUT}).
For any $u\in U \subseteq T, $ we clearly have  $ \|x- \Pi_ T(x) \|_2 \leq \|x-u\|_2  .$
Thus  by the  triangle inequality and (\ref{PPII}),  we have
\begin{eqnarray} \label{error-zhao}
\|x- \Pi_ \Omega(x) \|_2    & \leq & \| x- \Pi_ T(x) \|_2 + \| \Pi_ T(x) - \Pi_ \Omega(x)  \|_2  \nonumber \\
& \leq &   \| x- u \|_2  + \sqrt{ d^{\cal H} (\Omega ,
T)  \|x- \Pi_ \Omega (x) \|_2}
 \end{eqnarray}  for any $x \in \mathbb{R}^n$ and   $ u \in U. $
Note that the quadratic equation
$t^2  = \alpha + \sqrt{\beta} t $ in $t$, where $\alpha\geq 0 $ and $ \beta \geq 0$, admits a unique nonnegative root
$t^* = \frac{\sqrt{\beta} + \sqrt{\beta+4\alpha} }{2}.$ Thus, by setting $$t=\sqrt{\|x- \Pi_ \Omega (x) \|_2}, ~ \alpha =  \| x- u \|_2, ~  \beta = d^{\cal H} (\Omega ,
T), $$ it immediately follows from (\ref{error-zhao})   that
$$ \sqrt{\|x- \Pi_ \Omega(x) \|_2 }  \leq \frac{\sqrt{d^{\cal H } (\Omega, T)} + \sqrt{d^{\cal H} (\Omega,T)  +4\| x- u \|_2 } }{2},  $$ which implies that
$$
\|x- \Pi_ \Omega(x) \|_2 \leq  \left(\frac{\sqrt{d^{\cal H } (\Omega,T)} + \sqrt{d^{\cal H} (\Omega,T)  +4\| x- u \|_2 } }{2}\right)^2 \leq d^{\cal H}(\Omega,T)+ 2 \| x- u \|_2,     $$
  where the last inequality follows from the fact $(t_1+t_2)^2 \leq 2(t_1^2+t_2^2) $ for any numbers $t_1$ and $t_2.$  ~ $\Box$

  \vskip 0.07in


\vskip 0.07in

To show the next two technical results, let us first define a set which is a relaxation of the solution set $S^*$ of problem (\ref{L1}).  Note that $$S^*= \{x \in \mathbb{R}^n:  ~ \|x \|_1\leq \varrho^* ,  ~ u=M^T(Ax
-y)/\tau, ~ u\in \mathfrak{B}^\phi \} . $$ Replacing $\mathfrak{B}^\phi$ with the polytope $ {\cal Q}_{\epsilon_j}, $ which is an outer approximation of $\mathfrak{B}^\phi, $   yields the  set
\begin{equation} \label{S-epsilon}  S_{\epsilon_j}     =
\{x \in \mathbb{R}^n:  ~ \|x \|_1\leq \varrho^* ,  ~ u=M^T(Ax
-y)/\tau, ~ u\in {\cal Q}_{\epsilon_j} \}. \end{equation}
   Clearly,  $S^* \subseteq S_{\epsilon_j}$ due to the fact $\mathfrak{ B}^\phi \subseteq {\cal Q}_{\epsilon_j} . $

 \vskip 0.07in

\textbf{Lemma 4.2.} \emph{Let $ S_{\epsilon_j} $ be the set defined in (\ref{S-epsilon}) where $ {\cal Q}_{\epsilon_j} $ is given in (\ref{PJ-01}).
For every $j,$  let $x_{\epsilon_j} $ be an arbitrary point in $S_{\epsilon_j}. $ Then every accumulation point $\widehat{x} $ of the sequence $ \{x_{\epsilon_j} \}_{j\geq 1}  $ satisfies  $ \phi(M^T (A\widehat{x}-y)) \leq \tau.$  }

\vskip 0.07in

\emph{Proof.}    Recall that ${\cal Q}_{\epsilon_j}$  is  represented as (\ref{QQQ}), i.e.,  $$ {\cal Q}_{\epsilon_j} = \{u\in \mathbb{R}^q:  (a^i)^T   u \leq 1 \textrm{ for all } a^i \in \mathcal{Y}(\Gamma_{{\cal Q}_{\epsilon_j}}) \} .$$  Note that  $ x_{\epsilon_j}  \in S_{\epsilon_j}$
for any $j\geq  1 . $  Then for any $j \geq 1,$ we see from (\ref{S-epsilon}) that
\begin{equation} \label{ggaa01} \left\|x_{\epsilon_j} \right\|_1\leq \varrho^* ,  ~ (a^i)^T [M^T(A x_{\epsilon_j}- y )] \leq \tau   \textrm{ for all } a^i \in \mathcal{Y}(\Gamma_{{\cal Q}_{\epsilon_j}}). \end{equation}
  The   sequence $ \{x_{\epsilon_j} \} _{j\geq 1} $ is bounded. Let $\widehat{x} $ be any accumulation point of the sequence $ \{x_{\epsilon_j} \} _{j\geq 1}. $   Passing through to a subsequence if necessary, we may assume that  $x_{\epsilon_j} \to \widehat{x} $ as $j\to \infty. $  We prove the lemma by contradiction. Assume that $\phi(M^T(A \widehat{x}   -y)) > \tau.$
Then we define
  $$\widehat{\sigma}  :=  \frac{  \phi(M^T(A \widehat{x}   -y)) - \tau  }{\tau}, $$ which  is a positive constant under the assumption.     Since $\epsilon_j \to 0$ as $j\to \infty,$  there exists an integer number $j_0$ such that  $\epsilon_j < \widehat{\sigma}  $  for any $j\geq j_0.$
 By the definition of $\Gamma_ {{\cal Q}_{\epsilon_j}},$ we see that
    \begin{equation} \label{cont}       \mathcal{Y}  (\Gamma_ {{\cal Q}_{\epsilon_{j'}}}   ) \subseteq   \mathcal{Y}  (\Gamma_ {{\cal Q}_{\epsilon_j}}   )\textrm{  for any  }  j\geq j'\geq  j_0.   \end{equation}
     Thus   for
any  fixed integer number $j' \geq j_0,$ the following holds for all $j
\geq  j' :$
  $$ \sup_{a^i \in \mathcal{Y}(\Gamma_ {{\cal Q}_{\epsilon_{j'}}})  }  (a^i)^T [M^T(Ax_{\epsilon_j}-y)]    \leq  \sup_{a^i \in \mathcal{Y}(\Gamma_ {{\cal Q}_{\epsilon_j}}  ) }   (a^i)^T [M^T (Ax_{\epsilon_j}-y)]   \leq \tau ,
$$ where the first inequality follows from (\ref{cont}) and the second inequality follows from (\ref{ggaa01}).  For every fixed $j' \geq j_0$, noting that  $x_{\epsilon_j}\to \widehat{x}  $ as $j\to \infty,$  it follows from the  above inequality  that
\begin{equation}\label{J4}  \sup_{a^i \in  \mathcal{Y}(\Gamma_ {{\cal Q}_{\epsilon_{j'}} }) } (a^i)^T [M^T(A\widehat{x}
-y)] \leq \tau.  \end{equation}
 Consider the vector   $\widehat{a} = M^T(A\widehat{x}   -y)/\phi(M^T(A \widehat{x}  -y))  $ which is on the surface of the unit ball $ \mathfrak{B}^\phi.$
 Note that $ \epsilon_{j'}
<  \widehat{\sigma}.$ Applying Lemma 2.3 to $ {\cal Q}_{\epsilon_{j'}},$
we conclude that for the vector $\widehat{a} ,$  there is a vector $a^i \in
\mathcal{Y}(\Gamma_ {  {\cal Q}_{\epsilon_{j'} }}) $    such that
$   (a^i)^T \widehat{a}  \geq  \frac{1}{1+\epsilon_{j'}} > \frac{1}{1+\widehat{\sigma}}, $    which implies   that
\begin{equation}  \label{TTSS} (a^i)^T [M^T(A\widehat{x} -y) ] >\frac{\phi(M^T(A\widehat{x}  -y ))}{1+\widehat{\sigma}}   = \tau,\end{equation}  where the equality follows from the definition of $\widehat{\sigma}.$   Since $a^i \in \mathcal{Y}(\Gamma_ {  {\cal Q}_{\epsilon_{j'} }}), $ the inequality (\ref{TTSS}) contradicts (\ref{J4}). Therefore, for any accumulation point $\widehat{x}$ of the sequence $\{x_{\epsilon_j}\}_{j\geq 1}, $  we must have that $ \phi(M^T(A \widehat{x}  -y)) \leq \tau. $   ~  $ \Box $

\vskip 0.07in

\textbf{Lemma 4.3.} \emph{Let $S^*$ be the solution set of problem (\ref{L1}) and let $ S_{\epsilon_j} $ be the set given in (\ref{S-epsilon}).
 Then
$d^{\cal H} (S^*, S_{\epsilon_j}) \to 0 \textrm{
as } j\to \infty. $}

\vskip 0.07in

\emph{Proof.}
 Since  $S^* \subseteq S_{\epsilon_j},$
by the definition of Hausdorff metric, we see that \begin{equation}
\label{pipi01}  d^{\cal H} (S^*, S_{\epsilon_j}) =
\sup_{x\in S_{\epsilon_j} }
 \inf_{z\in S^*}\|x-z\|_2 = \sup_{x\in S_{\epsilon_j} }  \| x-\Pi_ {S^*} (x)\|_2. \end{equation}   Note that  $S^*$ and $ S_{\epsilon_j}$ are compact convex sets and  the projection operator $\Pi_ {S^*} (x) $ is
continuous in $ \mathbb{R}^n. $   For every polytope $S_{\epsilon_j},$ the
 superimum in (\ref{pipi01})  can be  attained. Thus there
exists a point in $ S_{\epsilon_j},$ denoted by $x_{ \epsilon_j},$ such that
 $ d^{\cal H} (S^*, S_{\epsilon_j}) = \left\|x_{\epsilon_j}  - \Pi_ {S^*}
 ( x_{\epsilon_j}) \right\|_2.
$
Note that  $ S^* \subseteq S_{ \epsilon_{j+1} }
\subseteq S_{\epsilon_j} $ for any $j\geq 1.$   The sequence
 $ \{d^{\cal H} (S^*, S_{\epsilon_j})
\}_{j \geq 1} $  is   non-increasing and nonnegative. The limit $ \lim_{j\to \infty} d^{\cal H} (S^*, S_{\epsilon_j})
   $ exists.
Passing through to subsequence if necessary, we may assume that the sequence $\{x_{\epsilon_j} \}_{j\geq 1} $ tends to $\widehat{x} . $ Note that $x_{\epsilon_j} \in S_{\epsilon_j} $ which indicates that $ \|x_{\epsilon_j} \|_1 \leq \varrho^* $ and hence $\|\widehat{x}\|_1 \leq \varrho^* . $  By Lemma 4.2,    $ \widehat{x}  $  must satisfy that  $ \phi(M^T(A \widehat{x} -y)) \leq \tau  $ which, together with   $\|\widehat{x} \|_1\leq \varrho^* , $
  implies that  $\widehat{x}  \in S^*.$ As a result, $ \Pi_ {S^*} (\widehat{x} )=\widehat{x} .$  Therefore,
  $$   \lim_{j\to \infty} d^{\cal H} (S^*, S_{\epsilon_j}) = \lim_{j\to \infty}
 \| x_{\epsilon_j}- \Pi_ {S^*}(
   x_{\epsilon_j})\|_2  =\| \widehat{x} -\Pi_ {S^*} ( \widehat{x} )\|_2 =0 ,  $$  as desired.  ~~ $
\Box $

\vskip 0.07in

Throughout the remainder of this section, let $   \delta  $ be any  fixed small constant (e.g.,  a sufficiently small constant in $(0, \tau) $). By Lemma 4.3. there is a $j_0$ such that  $S_{\epsilon_{j_0}}$ defined in (\ref{S-epsilon}) achieves
   \begin{equation} \label{DLT} d^{\cal H} (S^*, S_{\epsilon_{j_0}}) \leq \delta .  \end{equation}

  In the reminder of this section,   we focus on the  fixed polytope   ${\cal Q}_{\epsilon_{j_0}},$ as an approximation of $\mathfrak{B}^\phi. $ We use
$\widehat{n}$ to denote the number of the columns of
$\Gamma_{{\cal Q}_{\epsilon_{j_0}}}$ and use  $\widehat{\mathbf{e}}$ to
denote the vector of ones in $\mathbb{R}^{\widehat{n}} $ to distinguish the vector  of ones in other spaces.
Replacing $\mathfrak{B}^\phi $ in  (\ref{l1bb}) with $ {\cal Q}_{\epsilon_{j_0}} $
leads to the following relaxation of problem (\ref{L1}):
\begin{eqnarray} \label {PPJJa} \varrho^* _{j_0}: & = & \min_{(x,u)}\{\|x\|_1: ~ u=M^T(Ax-y)/\tau, ~ u\in
{\cal Q}_{\epsilon_{j_0}}\}  \nonumber \\
& =  & \min_x\{\|x\|_1:
~ (\Gamma_{ {\cal Q}_{\epsilon_{j_0}} } )^T [M^T(Ax-y) ] \leq \tau
\widehat{\mathbf{e}}\} , \end{eqnarray} where $\varrho^* _{j_0} $ denotes the optimal value of the above optimization problem.  Clearly, $\varrho^* _{j_0}  \leq \varrho^*    $ due to the
fact $\mathfrak{B}^\phi \subseteq {\cal Q}_{\epsilon_{j_0}}. $ Let $S^*_{\epsilon_{j_0}}$ denote the set of optimal solutions of (\ref{PPJJa}), i.e.,  $$
S^*_{\epsilon_{j_0}}  =  \{ x\in \mathbb{R}^n : ~
\|x\|_1 \leq \varrho^* _{j_0}, ~
u=M^T (Ax-y)/\tau,  ~ u\in
 {\cal Q}_{\epsilon_{j_0}} \}. $$  By (\ref{S-epsilon}), we immediately see
that $S^*_{\epsilon_{j_0}} \subseteq  S_{\epsilon_{j_0}}  $ since $\varrho^* _{j_0}  \leq \varrho^* .  $
 The  problem (\ref{PPJJa})   can be written
as
$$  \min_{(x,t)} \left\{ \widetilde{\textbf{e}}^T t: ~ x\leq t, ~ -x\leq t, ~ t\geq 0, ~(\Gamma_{{\cal Q}_{\epsilon_{j_0}}})^T [M^T(Ax-y)] \leq
\tau \widehat{\mathbf{e}} \right\}, $$ where $\widetilde{\mathbf{e}}$ is the vector of ones in $ \mathbb{R}^n.$ Clearly, $t=|x|$ at any optimal solution of the problem.  The Lagrangian dual of the above problem  is given as
\begin{eqnarray} \label{ddll2} &  \max  &  - \left[\tau \widehat{\mathbf{e}} + (M \Gamma_{{\cal Q}_{\epsilon_{j_0}}})^Ty \right]^T w_3   \\
& \textrm{s.t.} &  (A^T M \Gamma_{{\cal Q}_{\epsilon_{j_0}}})  w_3
+w_1-w_2=0, ~ w_1+w_2 \leq \widetilde{\mathbf{e}}, ~w_1, ~ w_2, ~ w_3\geq 0.  \nonumber
\end{eqnarray}
 By the KKT optimality
condition, the solution set of (\ref{PPJJa}) can be completely characterized.

\vskip 0.07in

\textbf{Lemma 4.4.} \emph{$x^* \in \mathbb{R}^n$ is an optimal
solution of (\ref{PPJJa})  if and only if there exist vectors $t^*,
w_1^*, w_2^* \in \mathbb{R}^n_+ $ and $
 w_3^* \in \mathbb{R}^{\widehat{n}} _+  $   such that $ (x^*, t^*,
w_1^*, w_2^*, w_3^*) \in {\mathfrak{D}}^{(1)}, $ where ${\mathfrak{D}}^{(1)}$ is the set of vectors  $(x,t,  w_1,w_2, w_3)$ satisfying the following system:
 \begin{equation} \label{TTTa} \left\{\begin{array} {l}
   x\leq t, ~ -x\leq t, ~  (  \Gamma_{{\cal Q}_{\epsilon_{j_0}}})^T [M^T(Ax-y)]\leq \tau \widehat{\mathbf{e}},
        \\
    w_1+w_2\leq \widetilde{\mathbf{e}},
   ~   A^T  M  \Gamma_{{\cal Q}_{\epsilon_{j_0}}} w_3 + w_1-w_2 =0,  \\
    \widetilde{\mathbf{e}}^T t =- \left[\tau \widehat{\mathbf{e}}+ (  M \Gamma_{{\cal Q}_{\epsilon_{j_0}}})^T   y  \right] ^T w_3,
\\
    (t, ~ w_1, ~ w_2,  ~ w_3)  \geq 0.
\end{array}
\right.
\end{equation}
 }

\vskip 0.07in
By optimality, we see that  $t=|x| $ for any $ (x,t, w_1,w_2,  w_3) \in \mathfrak{D}^{(1)}. $
  We write (\ref{TTTa})  as
$$ {\mathfrak{D}}^{(1)}  = \{z= (x,t, w_1,w_2,  w_3):   ~ \overline{M}^1 z\leq \overline{b}^1,  ~ \overline{M}^2  z = \overline{b}^2 \},
$$  where $\overline{b}^2=0$, $\overline{b}^1 = [
              0,
              0,
              \widetilde{\mathbf{e}}^T,
             ((M \Gamma_{{\cal Q}_{\epsilon_{j_0}}})^Ty +\tau \widehat{\mathbf{e}})^T,
              0,
              0,
              0,
              0]^T$ and
 \begin{equation} \label {MB++}  \overline{M}^1 = \left[
            \begin{array}{ccccc}
              I & -I & 0 & 0 & 0 \\
              -I & -I & 0 & 0 & 0 \\
              0 & 0 & I & I & 0 \\
              (M \Gamma_{{\cal Q}_{\epsilon_{j_0}}})^T A & 0 & 0 & 0 & 0 \\
              0 & -I & 0 & 0 & 0 \\
              0 & 0 & -I & 0 & 0 \\
              0 & 0 & 0 & -I & 0 \\
              0 & 0 & 0 & 0 & -\widehat{I}  \\
            \end{array}
          \right]
 ,
  \end{equation}
\begin{equation} \label{M2++}  \overline{M}^2 =  \left[
                \begin{array}{ccccc}
                  0 &  0 & I  & -I  & A^T M \Gamma_{{\cal Q}_{\epsilon_{j_0}}} \\
                  0 & \widetilde{\mathbf{e}}^T & 0  & 0  &  \tau \widehat{\mathbf{e}} ^T +   y^T M \Gamma_{{\cal Q}_{\epsilon_{j_0}}} \\
                \end{array}
              \right],
\end{equation}
where $I \in \mathbb{R}^{n\times n} $ and $\widehat{I}  \in \mathbb{R}^{\widehat{n}\times \widehat{n}} $  are identity
matrices  and $0$'s are zero-matrices with suitable sizes.
The main   result in this section is given as follows.

 \vskip 0.07in

\textbf{Theorem 4.5.}   \emph{
Given the problem data $(A, M, y, \tau),$ where  $ A \in \mathbb{R}^{m\times n} $ ($m<n $) and $ M \in \mathbb{R}^{m\times q} $ ($m\leq q$)  with
 $ \textrm{rank}(A)= \textrm{rank}(M)=m. $ Let $\delta  \in (0, \tau) $ be a  fixed constant and let $ {\cal Q}_{\epsilon_{j_0}}$ be the fixed polytope   such that (\ref{DLT}) is achieved.  Suppose that $A^T$ satisfies the weak RSP of order
$k. $
   Then for any $x \in \mathbb{R}^n, $
    there is a solution $x^*$ of (\ref{L1})
     approximating  $x$ with error \begin{equation}\label{ll2-error} \|x-x^*\|_2 \leq  \delta +
       2\overline{\gamma}
\left\{ \widehat{n}\left(\phi(M^T(Ax-y))-\tau\right)^+   +2\sigma_k(x)_1 + c
\tau  + c  \phi(M^T(Ax-y)) \right\} , \end{equation}
where $ c
 $ is the constant given in Theorem 3.2,  and $\overline{\gamma}  = \sigma_{\infty, 2} (\overline{M}^1, \overline{M}^2) $ is
 the Robinson
  constant determined by  $(\overline{M}^1, \overline{M}^2)$   in (\ref{MB++}) and  (\ref{M2++}).
 Moreover,  for any $x$ satisfying the constraint of (\ref{L1}),
     there is an optimal solution $x^*$ of (\ref{L1})
     approximating $x$ with error  }
     \begin{equation} \label{error-45} \|x-x^*\|_2 \leq  \delta + 2\overline{\gamma}  \left\{ 2\sigma_k(x)_1 +
 c   \tau + c \phi(M^T(Ax-y))   \right\}  \leq  \delta + 4\overline{\gamma}  \left\{  \sigma_k(x)_1 +
 c   \tau    \right\} . \end{equation}

\emph{Proof.}  Let $x$ be any given vector in $\mathbb{R}^n. $   Let $S$ be the support set of the $k$-largest
entries of $|x|.$  Let $S' =\{ i\in S: x_i> 0\} $
 and $ S'' = \{i\in S: x_i <  0\} .$ Clearly, $|S'|+|S''|\leq |S| \leq k,$ and $S'$ and $S''$ are disjoint.
 Since $A^T$ satisfies the weak RSP
 of order $k$, there exists a vector  $ \zeta= A^Tu^*  $ for some
 $u^*\in \mathbb{R}^m $   satisfying
 $ \zeta_i= 1 \textrm{ for }i\in S',  ~ \zeta_i= -1 \textrm{ for  }i\in S'', ~ \textrm{ and }
 |\zeta_i| \leq 1\textrm{ for } i \notin  S'\cup S'' . $
For the  fixed constant $\delta \in (0, \tau),$  there exists an integer number $j_0   $
such that the polytope ${\cal Q}_{\epsilon_{j_0}}$, represented as (\ref{QQQ}),
  ensures that the set $ S_{\epsilon_{j_0}}, $ defined  by (\ref{S-epsilon}),  achieves  the  bound (\ref{DLT}).    We now construct a
feasible solution $(\widetilde{w}_1, \widetilde{w}_2,
\widetilde{w}_3)$ to the problem (\ref{ddll2}). Set
$(\widetilde{w}_1)_i= 1\textrm{ and } (\widetilde{w}_2)_i=0\textrm{
for all }i\in S', ~ (\widetilde{w}_1)_i=0\textrm{ and
}(\widetilde{w}_2)_i= 1\textrm{ for all }i\in S'' , $ and  $
(\widetilde{w}_1)_i =(|\zeta_i|+\zeta_i)/2 $ and $  (\widetilde{w}_2)_i
=(|\zeta_i|-\zeta_i)/2 \textrm{ for all }i\notin S'\cup S''. $ This
 choice  ensures that $ (\widetilde{w}_1, \widetilde{w}_2)\geq
0,  ~  \widetilde{w}_1 + \widetilde{w}_2 \leq \widetilde{\textbf{e}} $ and $
\widetilde{w}_1 -\widetilde{w}_2  =\zeta. $

We now construct the
vector $\widetilde{w}_3. $
    By the definition of ${\cal Q}_{\epsilon_{j_0}}$, we see that
    $  \{\pm \textbf{e}_i: ~i=1, \dots,  q\} \subseteq \mathcal{Y}  (\Gamma_{ {\cal
    Q}_{\epsilon_{j_0} }}),
    $ the set of column vectors of  $\Gamma_{{\cal
    Q}_{\epsilon_{j_0} }}$ with cardinality $| \mathcal{Y}  (\Gamma_{{\cal
    Q}_{\epsilon_{j_0} }})| = \widehat{n}.$
  It is not difficult to show that
    there exists a
vector $\widetilde{w}_3\in \mathbb{R}^{\widehat{n}}_{+}$ satisfying
$M \Gamma_{{\cal Q}_{\epsilon_{j_0}} } \widetilde{w}_3= - u^* $.
First, since $M$ is a full row rank matrix, there exists an $m\times m$ invertible submatrix $M_\mathfrak{J}$ with $|\mathfrak{J}|=m $ consisting of $m$ independent columns in $M$. Then by choosing $\widetilde{h}\in \mathbb{R}^q$ such that $\widetilde{h}_i=0$ for all $ i\notin  \mathfrak{J}$ and $ \widetilde{h}_\mathfrak{J} =-M_{\mathfrak{J}}^{-1}u^* $ which implies that $M\widetilde{h} =-u^*.$ We now find $\widetilde{w}_3$ such that $\Gamma_{{\cal Q}_{\epsilon_{j_0}} } \widetilde{w}_3 =\widetilde{h}. $
In fact, without
loss of generality, we assume that $ \{-\textbf{e}_i: ~i=1, \dots, q\}$
are arranged as the first $q$ columns of
$\Gamma_{{\cal Q}_{\epsilon_{j_0}} }  $ and  $ \{\textbf{e}_i: ~i=1,
\dots,  q\}$ are arranged as the second $q$ columns of
$\Gamma_{{\cal Q}_{\epsilon_{j_0}} }. $ For every $i=1, \dots , q,$ if
$\widetilde{h}_i\geq 0,$ then we set  $ (\widetilde{w}_3)_i= \widetilde{h}_i; $ otherwise,
if $\widetilde{h}_i< 0,$ then we set   $ (\widetilde{w}_3)_{q+i}=  - \widetilde{h}_i.  $ All
remaining entries of $ \widetilde{w}_3  \in \mathbb{R}^{\widehat{n}}
$ are set to be zero. By this choice of $ \widetilde{w}_3,$ we see
that $\widetilde{w}_3 \geq 0$ satisfying that $\Gamma_{{\cal Q}_{\epsilon_{j_0}} }
\widetilde{w}_3= - \widetilde{h}  $ and
\begin{eqnarray}  \|\widetilde{w}_3\|_1  & =  & \|\widetilde{h}\|_1 =\|\widetilde{h}_\mathfrak{J}\|_1 =   \|M_{\mathfrak{J}}^{-1} u^* \|_1  = \|M_{\mathfrak{J}}^{-1} (AA^T)^{-1}A \zeta\|_1   \nonumber \\  & \leq &
\|M_{\mathfrak{J}}^{-1}(AA^T)^{-1}A \|_{\infty\to 1} \|\zeta\|_\infty \leq  c ,   \label{EEE1}
\end{eqnarray}
where $  c   $ is the constant given in Theorem 3.2.  By the triangle inequality and the fact $\phi^* (\mathbf{e}_i)=1, i=1, \dots, q,$    we have
\begin{equation}\label{6666}
\phi^* (\widetilde{h})= \phi^* (\sum_{j\in \mathfrak{J}} \widetilde{h}_j \mathbf{e}_j) \leq  \sum_{j\in \mathfrak{J}} \phi^* (\widetilde{h}_j \mathbf{e}_j) = \sum_{j\in \mathfrak{J}} |\widetilde{h}_j | \phi^* ( \mathbf{e}_j) = \sum_{j\in \mathfrak{J}} |\widetilde{h}_j | =\|\widetilde{h}\|_1 \leq c,
\end{equation}
where the last inequality follows from (\ref{EEE1}).

 For the vector $(x,t, \widetilde{w}_1, \widetilde{w}_2, \widetilde{w}_3) $ with $t=|x|,$  applying Lemma 2.1 with $(M^1, M^2)= (\overline{M}^1, \overline{M}^2) ,$ where $\overline{M}^1$ and $ \overline{M}^2$ are given in  (\ref{MB++}) and (\ref{M2++}),
  yields  a point in $ \mathfrak{D}^{(1)}  ,$  denoted by
$(\widehat{x}, \widehat{t}, \widehat{w}_1, \widehat{w}_2,
\widehat{w}_3), $ such that
\begin{equation} \label {hoff-44} \left\|
   \left[
      \begin{array}{c}
        x \\
        t \\
        \widetilde{w}_1 \\
        \widetilde{w}_2 \\
        \widetilde{w}_3 \\
      \end{array}\right] - \left[
       \begin{array}{c}
         \widehat{x} \\
        \widehat{t} \\
         \widehat{w}_1 \\
         \widehat{w}_2 \\
         \widehat{w}_3 \\
       \end{array}
     \right]
 \right\|_2 \leq  \overline{\gamma}  \left\| \left[\begin{array}{c}
(x-t)^+ \\
   (-x- t)^+ \\
   \left(( \Gamma_{{\cal Q}_{\epsilon_{j_0}} })^T[M^T(Ax-y)]-\tau \widehat{\mathbf{e}} \right)^+ \\
    ( \widetilde{w}_1+\widetilde{w}_2-\widetilde{\mathbf{e}})^+ \\
    A^T M \Gamma_{{\cal Q}_{\epsilon_{j_0}} } \widetilde{w}_3+ \widetilde{w}_1-\widetilde{w}_2  \\
     \widetilde{\mathbf{e}}^T t + \left(\tau \widehat{\mathbf{e}} +  ( M \Gamma_{{\cal Q}_{\epsilon_{j_0}} })^T y \right) ^T\widetilde{w}_3 \\
        (V)^+
            \end{array}
             \right] \right\|_1. \end{equation}
    where  $(V)^+=
( ( -t)^+, (-\widetilde{w}_1)^+,
          (-\widetilde{w}_2)^+,
          (-\widetilde{w}_3)^+),
      $ and  $\overline{\gamma}  =\sigma_{\infty, 2} (\overline{M}^1, \overline{M}^2) $ is the Robinson's constant
    determined by  $(\overline{M}^1, \overline{M}^2)$  in (\ref{MB++}) and (\ref{M2++}). By the nonnegativity of $ (t, \widetilde{w}_1,
\widetilde{w}_2, \widetilde{w}_3), $ we   have that  $(V)^+ =0.$
Since $t=|x|, $ we have  $(x-t)^+ =
   (-x- t)^+   =0  . $ Since $ (\widetilde{w}_1,
\widetilde{w}_2, \widetilde{w}_3) $ satisfies the constraints of (\ref{ddll2}), we
have      $( \widetilde{w}_1+\widetilde{w}_2-\widetilde{\mathbf{e}})^+ =0$ and $A^T M  \Gamma_{{\cal Q}_{\epsilon_{j_0}} } \widetilde{w}_3+ \widetilde{w}_1-\widetilde{w}_2 =0.$
Thus the inequality  (\ref{hoff-44}) is reduced to
\begin{equation}  \label {error-aa} \|x-\widehat{x} \|_2 \leq \overline{\gamma}  \left\{ \left\|
\left((  \Gamma_{{\cal Q}_{\epsilon_{j_0}} })^T[M^T(Ax-y)]-\tau \widehat{\mathbf{e}}  \right)^+ \right\|_1+
\left|\widetilde{\mathbf{e}}^T t +  \left[\tau \widehat{\mathbf{e}} + ( M \Gamma_{{\cal Q}_{\epsilon_{j_0}}
})^T y  \right]^T \widetilde{w}_3\right|\right\}.
\end{equation}
 Recall
that $\phi^*(a^i) =1 $ for every $a^i \in \mathcal{Y}(
\Gamma_{{\cal Q}_{\epsilon_{j_0}} }). $ Thus $$ (a^i)^T (M^T(Ax-y)) \leq \phi^*(a^i)
\phi(M^T(Ax-y)) = \phi(M^T(Ax-y)) .$$ This implies that $ \left((a^i)^T (M^T(Ax-y))- \tau \right)^+ \leq
\left(\phi(M^T(Ax-y))-\tau \right)^+ , $ and hence
$$ \left((
\Gamma_{{\cal Q}_{\epsilon_{j_0}} })^T (M^T (Ax-y)) -\tau \widehat{\mathbf{e}} \right)^+ \leq (\phi(M^T(Ax-y))
-\tau)^+ \widehat{\mathbf{e}}. $$ Therefore,
\begin{equation} \label{NNNN}  \left\|\left((\Gamma_{{\cal Q}_{\epsilon_{j_0}} })^T(M^T(Ax-y))-\tau \widehat{\mathbf{e}} \right)^+\right\|_1 \leq \widehat{n}
(\phi(M^T(Ax-y))-\tau )^+. \end{equation}   Note that  $x^T A^T u^*=  x^T \zeta  = \|x_S\|_1+ x_{\overline{S}}^T
\zeta_{\overline{S}}   $ and  $|x_{\overline{S}}^T
\zeta_{\overline{S}}| \leq \|x_{\overline{S}}\|_1
\|\zeta_{\overline{S}}\|_\infty \leq \|x_{\overline{S}}\|_1. $ Thus
$$ \left|\widetilde{\mathbf{e}}^T |x| - x^T A^T u^* \right|= \left| \widetilde{\mathbf{e}}^T |x| -\|x_S\|_1 - x_{\overline{S}}^T
\zeta_{\overline{S}} \right| \leq 2
\|x_{\overline{S}}\|_1 =2\sigma_k(x)_1.$$
Denote by $\vartheta= M^T(Ax-y) $ and note that  $ M \Gamma_{{\cal Q}_{\epsilon_{j_0}} }\widetilde{w}_3=-u^*  $ and $ \Gamma_{{\cal Q}_{\epsilon_{j_0}} }\widetilde{w}_3 =-\widetilde{h}. $   We have
 \begin{eqnarray} \label {error-bb}  |\widetilde{\mathbf{e}}^T t +[\tau \widehat{\mathbf{e}} + ( M \Gamma_{{\cal Q}_{\epsilon_{j_0}} })^T y
]^T \widetilde{w}_3 |
& = & |\widetilde{\mathbf{e}}^T |x|  + \tau  \widehat{e} ^T \widetilde{w}_3 + (x^T A^T  M
 - \vartheta^T)  \Gamma_{{\cal Q}_{\epsilon_{j_0}} } \widetilde{w}_3  |  \nonumber \\
& = & |\widetilde{\mathbf{e}}^T |x| - x^T A^T u^* + \vartheta^T \widetilde{h}  + \tau \widehat{\mathbf{e}}^T \widetilde{w}_3 | \nonumber \\
& \leq  & 2\sigma_k(x)_1  + |\vartheta^T
\widetilde{h}| + |\tau \widehat{\mathbf{e}}^T \widetilde{w}_3|   \nonumber \\
& \leq & 2\sigma_k(x)_1  +
 \phi(\vartheta) \phi^*(\widetilde{h})   + \tau  \|\widetilde{w}_3\|_1 \nonumber \\
& \leq & 2\sigma_k(x)_1  + c  \phi(M^T(Ax-y))+  c\tau ,
\end{eqnarray}
where the last inequality follows from (\ref{EEE1}) and (\ref{6666}).
Combination of (\ref{error-aa}), (\ref{NNNN}) and  (\ref{error-bb}) gives rise to
\begin{equation}\label{error-CC} \|x-\widehat{x}\|_2 \leq
     \overline{\gamma}
\left\{  {\widehat{n}} \left(\phi(M^T(Ax-y))-\tau\right)^+   +2\sigma_k(x)_1 + c
\tau  + c  \phi(M^T(Ax-y))  \right\}. \end{equation}

We now consider the three bounded convex sets   $S^*, S^*_{\epsilon_{j_0}} $ and $ S_{\epsilon_{j_0}} . $ By their definitions,  $ S^* \subseteq  S_{\epsilon_{j_0}} $ and $ S^*_{\epsilon_{j_0}}  \subseteq  S_{\epsilon_{j_0}} .  $
 Let
$x^* = \Pi_ {S^*} (x)  $  and $ \overline{x}  =\Pi_ {S_{\epsilon_{j_0}} } (x) .$  Note that  $\widehat{x} \in S^*_{\epsilon_{j_0}}.$ Applying Lemma 4.1 by setting $S=S^*,   U= S^*_{\epsilon_{j_0}}$ and $ T= S_{\epsilon_{j_0}},$  we conclude that
$$ \|x-x^*\|_2 \leq d^{\cal H} (S^*, S_{\epsilon_{j_0}}) + 2 \|x-\widehat{x}\|_2 \leq \delta+ 2 \|x-\widehat{x}\|_2 .$$
Combining  this inequality with (\ref{error-CC}) yields (\ref{ll2-error}), i.e.,
$$ \|x-x^*\|_2 \leq  \delta + 2\overline{\gamma}  \left\{\widehat{n} \left(\phi(M^T(Ax-y) )-\tau  \right)^+ +  2\sigma_k (x)_1 +
c \tau + c \phi(M^T(Ax-y))   \right\}. $$   When $x $ satisfies
$\phi(M^T(Ax-y)) \leq \tau $, the inequality above  is reduced to (\ref{error-45}). ~~  $ \Box $

\vskip 0.07in
Since every condition listed in Corollary 3.5  implies the weak RSP of order $k,$  we  immediately have the following result   for  DS  with a nonlinear constraint.

\vskip 0.07in

\textbf{ Corollary 4.6.} \emph{Let $A$ and $M$ be given as in Theorem 4.5. Let $\delta \in (0, \tau) $ be  a fixed constant and let $ {\mathcal Q}_{\epsilon_{j_0}}$ be the fixed polytope represented as (\ref{QQQ}) such that (\ref{DLT}) is achieved. Suppose that one of the following conditions holds: }
  (a)  \emph{$A$ (with $\ell_2$-normalized columns) satisfies  the mutual coherence property
$\mu_1(k) + \mu_1(k-1) <1; $}   (b)  \emph{RIP of order $2k$ with constant $\delta_{2k} < 1/\sqrt{2};  $ } (c) \emph{stable NSP of order $k$ with constant $\rho \in (0,1); $}
(d) \emph{robust NSP of order $k$  with  $  \rho' \in (0,1) $
and $\rho'' >0; $}
(e)  \emph{NSP of order $k; $}
(f) \emph{RSP of order $k $ of $A^T.$  Then the conclusions of Theorem 4.5 are valid for the DS problem (\ref{L1}).}

\section{The LASSO problem}

  In this section, we consider the    nonlinear  minimization problem  (\ref{GLASSO}) which is still called a  LASSO problem in this paper since it includes the standard LASSO   as a special case.
Let $\rho^*$ denote the optimal value of   (\ref{GLASSO}), i.e.,  \begin{equation}\label{Nonlinear-Lasso} \rho^*  = \min_x \{\phi(M^T(Ax-y)) :  ~ \|x\|_1  \leq \mu
\},
\end{equation} where the problem data $(M,A, y, \mu)$ is   given, and $ \phi(\cdot) $ is any norm with $\phi (\textbf{e}_i)=1$ and  $\phi^*(\textbf{e}_i)=1$ for $ i=1,\dots, q. $    In this section, we  focus on the nonlinear norm  $\phi$ in the sense that the inequality $\phi(x) \leq t$ cannot be represented as a finite number of linear inequalities or equalities, for instance, when $\phi $ is   the    $\ell_p$-norm with  $   p \in (1, \infty). $    We show that problem (\ref{GLASSO}) is also  stable under the weak RSP of order $k  $ of $A^T.$   As a result, the  stability theorem can  be established for LASSO with a broad range of matrix properties.
Problem (\ref{GLASSO}), i.e., (\ref{Nonlinear-Lasso}),  is equivalent to
\begin{equation}\label{l1-l2bb} \rho^*   = \min_{(x, \rho)} \{\rho : ~  \phi(M^T(Ax-y))\leq \rho, ~  \|x\|_1  \leq \mu
\}.
\end{equation}
 Let $\Lambda^*$   be the set of optimal
solutions of (\ref{GLASSO}),
which in terms of $\rho^* $ can be written as
$$ \Lambda^* =\{  x  \in \mathbb{R}^n  :  ~  \|x\|_1\leq \mu,  ~ \phi(M^T(Ax-y)) \leq \rho^*  \}. $$
Since the first  constraint in (\ref{l1-l2bb})  is nonlinear,   we   use the analytic method  in Section 4     to develop the stability result for problem  (\ref{GLASSO}).

Recall that ${\cal Q}_{\epsilon_j}$,   defined in (\ref{PJ-01}),  is a polytope approximation of   $\mathfrak{B}^\phi $ in $\mathbb{R}^q,$  and   is represented as $${\cal Q}_{\epsilon_j} =\{u\in \mathbb{R}^q:  (a^i)^T u\leq 1\textrm{ for all } a^i \in {\cal Y} (\Gamma_{{\cal Q}_{\epsilon_j}}) \}.$$ The   vectors $a^i $ of  $ {\cal Y} (\Gamma_{{\cal Q}_{\epsilon_j}}) $
are drawn on the surface of the dual unit ball, i.e.,  $a^i  \in \{a\in \mathbb{R}^q:  ~\phi^* (a)=1\}. $ Using this  approximation, we introduce a relaxation of the solution set $\Lambda^*: $
\begin{equation} \label{L-epsilon} \Lambda_{\epsilon_j}     =
\{ x \in \mathbb{R}^n :    ~  \|x \|_1\leq \mu ,  ~  (a^i)^T (M^T(Ax
-y))\leq   \rho^*\textrm{ for all }a_i \in  {\cal Y} (\Gamma_{{\cal Q}_{\epsilon_j}})  \}. \end{equation}
Clearly, $\Lambda^* \subseteq \Lambda _{\epsilon_j} $ for any $j\geq 1. $  Then we have the following lemma.
\vskip 0.07in

\textbf{Lemma 5.1.} \emph{Let $ \Lambda_{\epsilon_j} $ be   defined in (\ref{L-epsilon}). The following properties hold:
 }

(i) \emph{For every $j,$  let $ x_{\epsilon_j}  $ be an arbitrary point in $\Lambda_{\epsilon_j}. $ Then any accumulation point $  \widehat{x}  $ of the sequence $ \{x_{\epsilon_j}\}  $  satisfies that $\|\widehat{x} \|_1\leq \mu  $  and $ \phi(M^T(A\widehat{x}-y)) \leq  \rho^* ,$  i.e., $\widehat{x}\in \Lambda^*. $ }

(ii) $d^{\cal H} (\Lambda^*, \Lambda_{\epsilon_j}) \to 0 \textrm{
as } j\to \infty. $

\vskip 0.07in

\emph{Proof.} The proof is similar to that of  Lemmas 4.2 and 4.3.   Note
that $ x_{\epsilon_j}    \in \Lambda_{\epsilon_j}$
for any $j\geq  1.$ Thus for every $j,$ we have
\begin{equation} \label{ggaa*}     \left\|x_{\epsilon_j} \right\|_1\leq \mu,  ~ (a^i)^T (M^T(A x_{\epsilon_j}- y )) \leq \rho^*      ~~ \textrm{  for all }  a^i \in \mathcal{Y} (\Gamma_ {{\cal Q}_{\epsilon_j}})
 . \end{equation}
 Let $\widehat{x} $ be any accumulation point which clearly obeys $\|\widehat{x}\|_1 \leq \tau.$   Passing through to a subsequence if necessary, we may  assume that  $x_{\epsilon_j} \to \widehat{x} $ as $j\to \infty. $ Assume that $\phi(M^T( A \widehat{x}   -y )) > \rho^*.$   We now prove that this assumption leads to a contradiction.
  Under this assumption, we define
  $$ \sigma^*   :=   \left\{\begin{array} {ll} \frac{  \phi(M^T(A \widehat{x}   -y)) - \rho^*  }{\rho^*} &  \rho^* \not= 0, \\ 1 & \rho^*  =  0 ,  \end{array} \right.  $$ which is a positive constant.    By the definition of $\epsilon_j$, there exists an integer number $j_0$ such that $\epsilon_j <  \sigma^*   $  for any $j\geq j_0. $      By a  similar argument in the proof   of Lemma 4.3, for any $ j\geq j' \geq  j_0$, it follows from  (\ref{ggaa*}) and the fact $\mathcal{Y} (\Gamma_ {{\cal Q}_{\epsilon_{j'}}}) \subseteq \mathcal{Y} (\Gamma_ {{\cal Q}_{\epsilon_j}}) $   that
\begin{equation}\label{JJ''}  \sup_{a^i \in  \mathcal{Y} (\Gamma_ {{\cal Q}_{\epsilon_{j'}} }) } (a^i)^T  (M^T(A \widehat{x}
-y))  \leq \rho^*.   \end{equation}
Let $ \widehat{a}= M^T(A\widehat{x}  -y)/ \phi(M^T(A \widehat{x}  -y)),  $ which is on the surface of   $\mathfrak{B}^{\phi}.$
   Applying Lemma 2.3 to $ {\cal Q}_{\epsilon_{j'}}$ for $j'\geq  j_0,$
we see that for   $\widehat{a},$  there is a vector $a^i \in
\mathcal{Y}(\Gamma_ {  {\cal Q}_{\epsilon_{j'} }})  $ such that
$  (a^i)^T \widehat{a} \geq  \frac{1}{1+\epsilon_{j'}}  > \frac{1}{1+ \sigma^* }  , $  which  implies  that
$$ (a^i)^T  [M^T(A\widehat{x} -y)]   >   \frac{\phi(M^T( A\widehat{x}  -y ))}{1+ \sigma^* }    \geq   \rho^*,$$ where the second inequality  follows from the definition of $ \sigma^* . $
This contradicts (\ref{JJ''}). Therefore, we must have that $\phi(M^T(A\widehat{x}-y))\leq \rho^*.$ This, together with $\|\widehat{x}\|\leq \tau, $ implies that  $\widehat{x} \in \Lambda^*. $

  We now prove that $
d^{\cal H} (\Lambda^*, \Lambda_{\epsilon_j}) \to 0 $ as $ j \to
\infty. $ Since $\Lambda^* \subseteq \Lambda_{\epsilon_j}$,
 by the continuity of  $\Pi_ {\Lambda^*}  (\cdot) $ and compactness of  $\Lambda_{\epsilon_j}, $   there exists for each $\epsilon_j$ a point $\widetilde{x}_{\epsilon_j} \in \Lambda_{\epsilon_j} $ such that  \begin{equation}
\label{pipi}  d^{\cal H} (\Lambda^*, \Lambda_{\epsilon_j} ) =
\sup_{x\in \Lambda_{\epsilon_j} }
 \inf_{z\in \Lambda^*}\|x-z\|_2 = \sup_{x\in \Lambda_{\epsilon_j} }  \| x-\Pi_ {\Lambda^*} (x)\|_2 =  \left\|\widetilde{x}_{\epsilon_j}  - \Pi_ {\Lambda^*}
 ( \widetilde{x}_{\epsilon_j}) \right\|_2.  \end{equation}
We also note that  $ \Lambda^* \subseteq \Lambda_{\epsilon_{j+1}}
\subseteq \Lambda_{\epsilon_j} $ for any $j\geq 1.$   Thus
 $ \{d^{\cal H} (\Lambda^*, \Lambda_{\epsilon_j})
\}_{j\geq 1} $  is a non-increasing nonnegative sequence. Thus
$ \lim_{j\to \infty} d^{\cal H} (\Lambda^*, \Lambda_{\epsilon_j})  $ exists.   Since  the sequence $ \{\widetilde{x}_{\epsilon_j}  \} _{j\geq 1} $ is bounded,
 passing through to subsequence of $\{\widetilde{x}_{\epsilon_{j}}\} $ if necessary, we may assume that  $  \widetilde{x}_{\epsilon_{j}}  \to \widetilde{x}$ as $j\to \infty.$ By result (i),  $\widetilde{x}$ must be in $\Lambda^*.$  Therefore $ \Pi_ {\Lambda^*} (\widetilde{x} )=\widetilde{x} . $ It follows from (\ref{pipi}) that
  $   \lim_{j\to \infty} d^{\cal H} (\Lambda^*, \Lambda_{\epsilon_j})  =\| \widetilde{x}-\Pi_ {\Lambda^*} ( \widetilde{x} )\|_2 =0 .  $    ~~ $
\Box $

\vskip 0.07in

In the reminder of this section, let $  \delta $ be any fixed sufficiently small constant. By Lemma 5.1. there exists an integer number $j_0$ such that
 \begin{equation} \label{DLT1} d^{\cal H} (\Lambda^*, \Lambda_{\epsilon_{j_0}}) \leq \delta ,  \end{equation}
 where $ \Lambda_{\epsilon_{j_0}}$ is the set ({\ref{L-epsilon})
 determined by   ${\cal Q}_{\epsilon_{j_0}} .$  We use
$\widehat{n}=  |{\cal Y} (\Gamma_{{\cal Q}_{\epsilon_{j_0}}})| $ to denote the number of columns of
$\Gamma_{{\cal Q}_{\epsilon_{j_0}}}$ and    $\widehat{\mathbf{e}} $   the vector of ones in $\mathbb{R}^{\widehat{n}}. $   Thus ${\cal Q}_{\epsilon_{j_0}}$ is represented as   $$ {\cal Q}_{\epsilon_{j_0}} =\{u \in \mathbb{R}^q: ~ (\Gamma_{{\cal Q}_{\epsilon_{j_0}}})^T u\leq  \widehat{\mathbf{e}} \}. $$
We consider  the following relaxation of (\ref{l1-l2bb}):
\begin{equation} \label {PPJJ} \rho^*_{\epsilon_{j_0}}:  = \min_{(x,\rho)} \{ \rho: ~  \|x\|_1 \leq \mu,
~ (\Gamma_{{\cal Q}_{\epsilon_{j_0}}} )^T (M^T(Ax-y)) \leq \rho
 \widehat{\mathbf{e}}  \} , \end{equation} where $\rho^*_{\epsilon_{j_0}}  $ denotes the optimal value of the above  optimization problem.
 Clearly,  $
\rho^*_{\epsilon_{j_0}}  \leq \rho^*   $ due to the
fact that (\ref{PPJJ}) is a relaxation of (\ref{l1-l2bb}).
 Since  $\Gamma_{{\cal Q} _{\epsilon_{j_0}}}  $    includes $ \pm \textbf{e}_i, ~i=1,\dots , n$  as its columns, the variable $\rho$ in (\ref{PPJJ})  must be nonnegative. Let $$
\Lambda^*_{\epsilon_{j_0}}  =  \{ x\in \mathbb{R}^n : ~
\|x\|_1 \leq \mu, ~
  (\Gamma_{{\cal Q}_{\epsilon_{j_0}}} )^T (M^T(Ax-y)) \leq \rho^*_{\epsilon_{j_0}} \widehat{\mathbf{e}}  \} $$ be the set of optimal solutions of  (\ref{PPJJ}). Recall that
  $$\Lambda_{\epsilon_{j_0} }     =
\{ x \in \mathbb{R}^n :    ~  \|x \|_1\leq \mu ,  ~   (\Gamma_{{\cal Q}_{\epsilon_{j_0}}} )^T (M^T(Ax-y)) \leq \rho^* \widehat{\mathbf{e}}  \}.
  $$     Clearly,
  $\Lambda^*_{\epsilon_{j_0}} \subseteq \Lambda_{\epsilon_{j_0}}   $ due to fact  $\rho^*_{\epsilon_{j_0}} \leq \rho^*. $  The  problem  (\ref{PPJJ})  can be written
as
\begin{equation} \label{LLPP3}  \min_{(x,t, \rho)} \left\{ \rho:   ~ x\leq t, ~ -x\leq t,   ~ \widetilde{\mathbf{e}}^T t \leq \mu, ~(\Gamma_{{\cal Q}_{\epsilon_{j_0}}})^T [M^T(Ax-y)] \leq
\rho  \widehat{\mathbf{e}}  ,  ~(t, \rho) \geq 0 \right\}, \end{equation}
where $ \widetilde{\mathbf{e}} $ is still the vector of ones in $ \mathbb{R}^n. $ It is straightforward to  verify that the Lagrangian dual   of this   problem  is given as
\begin{eqnarray} \label{ddll} &  \max  &  - \mu w_3  - (y^T M \Gamma_{{\cal Q}_{\epsilon_{j_0}}}) w_4    \\
& \textrm{s.t.} &  A^T M \Gamma_{{\cal Q}_{\epsilon_{j_0}}} w_4
+w_1-w_2=0, ~ w_1+w_2 -w_3 \widetilde{\mathbf{e}} \leq 0, ~ \widehat{\mathbf{e}} ^T w_4 \leq 1, \nonumber \\
& &  w_1\in \mathbb{R}^n_+, ~ w_2 \in \mathbb{R}^n_+, ~ w_3\in \mathbb{R}_+, ~w_4\in \mathbb{R}^{\widehat{n}}_+.  \nonumber
\end{eqnarray}
  The next lemma  follows immediately from  the KKT optimality
condition of (\ref{LLPP3}) or (\ref{ddll}).

\vskip 0.07in

\textbf{Lemma 5.2.} \emph{$\overline{x}  \in \mathbb{R}^n$ is an optimal
solution of (\ref{PPJJ})  if and only if there exist vectors $\overline{t} ,
\overline{w}_1 , \overline{w}_2  \in \mathbb{R}^n_+,  $  $  \overline{\rho} \in \mathbb{R}_+, $ $  \overline{w}_3 \in \mathbb{R}_+ $ and $
 \overline{w}_4  \in \mathbb{R}^{\widehat{n}} _+  $   such that $ ( \overline{x}, \overline{t}, \overline{\rho},
\overline{w}_1,  \overline{w}_2, \overline{w}_3,   \overline{w}_4) \in \mathfrak{D}^{(2)},  $ where $ \mathfrak{D}^{(2)} $ is the set of  vectors  $ (x,t,\rho,  w_1,w_2, w_3, w_4)$ satisfying the following system:
 \begin{equation} \label{TTT3}  \left\{\begin{array} {l}
    x\leq t, ~ -x\leq t,    ~ \widetilde{\mathbf{e}}^T t \leq \mu, ~(\Gamma_{{\cal Q}_{\epsilon_{j_0}}})^T [M^T(Ax-y)] \leq
\rho  \widehat{\mathbf{e}} ,
      \\
     A^T M \Gamma_{{\cal Q}_{\epsilon_{j_0}}} w_4
+w_1-w_2=0, ~ w_1+w_2 -w_3 \widetilde{\mathbf{e}} \leq 0,     \\
   ~  \widehat{\mathbf{e}}  ^T w_4 \leq 1, ~ \rho =  - \mu w_3  - (y^T M  \Gamma_{{\cal Q}_{\epsilon_{j_0}}}) w_4    , \\
  ( t,\rho,  w_1,w_2, w_3, w_4)\geq 0.
\end{array}  \right.
 \end{equation}
  }

\vskip 0.07in
By optimality, it is evident that $t=|x| $ for any $ (x,t, \rho, w_1, w_2, w_3,  w_4) \in \mathfrak{D}^{(2)}. $
Clearly,  (\ref{TTT3}) can be written as
 \begin{equation} \label{TTT3a} \mathfrak{D}^{(2)}   = \{z= (x, t,\rho,  w_1,w_2, w_3, w_4):   ~ \widehat{M}^1 z\leq \widehat{b}^1,  ~ \widehat{M}^2  z = \widehat{b}^2 \},
 \end{equation}  where $\widehat{b}^2=0,$   $\widehat{b}^1 =  (0, 0 , \mu, y^T M \Gamma_{{\cal Q}_{\epsilon_{j_0}}}   , 0, 1, 0, 0, 0, 0, 0, 0)^T $  and
 \begin{equation} \label {M+b+}  \widehat{M}^1 = \left[
                                         \begin{array}{ccccccc}
                                           I & -I & 0 & 0 & 0 & 0 & 0 \\
                                          -I & -I & 0 & 0 & 0 & 0 & 0 \\
                                           0 & \widetilde{\mathbf{e}}^T & 0 & 0 & 0 & 0 & 0 \\
                                          (M \Gamma_{{\cal Q}_{\epsilon_{j_0}}})^T A & 0 & -\widehat{\mathbf{e}} &  0 & 0  & 0 & 0 \\
                                           0 & 0 & 0 & I & I & -  \widetilde{\mathbf{e}}  & 0 \\
                                           0 & 0 & 0 & 0 & 0 & 0 & \widehat{\mathbf{e}}^T \\
                                           0 & -I & 0 & 0 & 0 & 0 & 0 \\
                                           0 & 0 & -1 & 0 & 0 & 0 & 0 \\
                                           0 & 0 & 0 & -I & 0 & 0 & 0 \\
                                           0 & 0 & 0 & 0 & -I & 0 & 0 \\
                                           0 & 0 & 0 & 0 & 0 & -1 & 0 \\
                                           0 & 0 & 0 & 0 & 0 & 0 & -\widehat{I} \\
                                         \end{array}
                                       \right],
  \end{equation}
\begin{equation} \label{M++}  \widehat{M}^2 =  \left[
                \begin{array}{ccccccc}
                  0 &  0 & 0 & I  & -I  & 0 &  A^T  M \Gamma_{{\cal Q}_{\epsilon_{j_0}}} \\
                  0 & 0 & 1 & 0  & 0  & \mu &      y^T  M \Gamma_{{\cal Q}_{\epsilon_{j_0}}}
                \end{array}
              \right],
\end{equation}
where $I\in \mathbb{R}^{n\times n}$  and  $\widehat{I}\in  \mathbb{R}^{\widehat{n}\times \widehat{n}}$ are identity
matrices and $0$'s are zero matrices with suitable sizes.
We now prove the main result in this section.

 \vskip 0.07in

\textbf{Theorem 5.3.}    \emph{Let $\delta >0$ be any fixed sufficiently small constant, and let $ {\cal  Q}_{\epsilon_{j_0}}$ be the fixed polytope represented as (\ref{QQQ}) such that  (\ref{DLT1}) is achieved. Let the   data $(M, A, y,
\mu)$  in (\ref{GLASSO}) be given, where $\mu >0 , $
  $ A \in \mathbb{R}^{m\times n}$  ($m<n$)  and $M \in R^{m\times q} $ ($m\leq q$)   with $rank (A)= rank (M)=m.$  Suppose that  $A^T$ satisfies the weak RSP of order
$k. $
   Then for any $x \in \mathbb{R}^n, $
    there is a solution $x^*$ of (\ref{GLASSO})
     approximating  $x$ with error
     \begin{equation}\label{l2-error}    \|x-x^*\|_2 \leq  \delta + 2   \widehat{\gamma}
\left[ \left(\|x\|_1-\mu \right)^+  +  2 \phi (M^T(Ax-y) ) +  \frac{ |\mu -\|x\|_1 |+   2\sigma_k (x)_1 } {c }  \right],
  \end{equation}
where  $ c $ is the  constant given in Theorem 3.2,   and $ \widehat{\gamma}  = \sigma_{\infty, 2} (\widehat{M}^1, \widehat{M}^2) $ is
 the Robinson's constant determined by  $(\widehat{M}^1, \widehat{M}^2)$   in (\ref{M+b+}) and  (\ref{M++}).
 Moreover,  for any $x$ with $\|x\|_1\leq \tau, $
     there is an optimal solution $x^*$ of (\ref{GLASSO})
     approximating  $x$ with error \begin{equation} \label{1616}    \|x-x^*\|_2 \leq  \delta + 2   \widehat{\gamma}
\left[ 2 \phi (M^T(Ax-y))   +  \frac{|\mu - \|x\|_1|   +  2\sigma_k (x)_1 } {c   }  \right]    .
     \end{equation} }

\emph{Proof.}   Let $x$ be any vector in $\mathbb{R}^n .  $
Then set $t=|x|  $ and let $\rho = \phi(M^T(Ax-y))  $ which implies from  ${\cal Y} (\Gamma_{\mathcal{Q}_{\epsilon_{j_0}}}) \subseteq \{a\in \mathbb{R}^q: \phi^* (a)=1\} $ that \begin{equation} \label{5544} (\Gamma_{{\cal Q}_{\epsilon_{j_0}}} )^T (M^T (Ax-y)) \leq \rho \widehat{\mathbf{e}}. \end{equation}    Denote by $S$ the support set of the $k$-largest
entries of $|x|.$  Let $S' =\{ i\in S: x_i> 0\} $
 and $ S'' = \{i\in S: x_i <  0\} .$
 Since $A^T$ satisfies the weak RSP
 of order $k$, there exists a vector  $ \zeta= A^Tu^*  $ for some
 $u^*\in \mathbb{R}^m $   satisfying
 $ \zeta_i= 1 \textrm{ for }i\in S',  ~ \zeta_i= -1 \textrm{ for  }i\in S'', ~ \textrm{ and }
 |\zeta_i| \leq 1\textrm{ for } i \notin S'\cup S'' . $
 Let $c$ be the constant given in Theorem 3.2.
 We now construct a set of vectors
  $(\widetilde{w}_1, \widetilde{w}_2,
\widetilde{w}_3, \widetilde{w}_4)$ which satisfies the constraints (\ref{ddll}). First, we set $\widetilde{w}_3 = 1/c. $   Then we set
$$(\widetilde{w}_1)_i= 1/c \textrm{ and } (\widetilde{w}_2)_i=0\textrm{
for all }i\in S', $$ $$ (\widetilde{w}_1)_i=0\textrm{ and
}(\widetilde{w}_2)_i= 1/c \textrm{ for all }i\in S'' , $$   $$
(\widetilde{w}_1)_i =\frac{|\zeta_i|+\zeta_i}{2c} \textrm{  and   } (\widetilde{w}_2)_i
=\frac{|\zeta_i|-\zeta_i}{2c} \textrm{ for all }i\notin S'\cup S''. $$  This
 choice of $\widetilde{w}_1$ and $\widetilde{w}_2$ implies that $ (\widetilde{w}_1, \widetilde{w}_2)\geq
0,  ~  \widetilde{w}_1 + \widetilde{w}_2 \leq  \widetilde{w}_3 \widetilde{\textbf{e}}, $  and $
\widetilde{w}_1 -\widetilde{w}_2  =\zeta/c. $ We now construct the
vector $\widetilde{w}_4  $ as follows.
    By the definition of ${\cal Q}_{\epsilon_{j_0}}$, we see that
    $  \{\pm \textbf{e}_i: ~i=1, \dots,  q\} \subseteq \mathcal{Y}  (\Gamma_{{\cal Q}_{\epsilon_{j_0}}}).
    $
  Since $M$ has a full row rank matrix, there exists an $m\times m$  invertible square submatrix   $M_\mathfrak{J},$ where $  \mathfrak{J} \subseteq \{1, \dots, q\}$ with $|\mathfrak{J}|=m.$ We define the vector $\widetilde{g} \in \mathbb{R}^q $ as follows:   $(\widetilde{g})_\mathfrak{J}=M_{\mathfrak{J}}^{-1} u^* $ and $(\widetilde{g})_i=0$ for $i \notin \mathfrak{J}.$ Clearly, we have $ M \widetilde{g} =u^*.$     It is not difficult to show that
    there exists a
vector $\widetilde{w}_4\in \mathbb{R}^{\widehat{n}}_{+}$ satisfying
$\Gamma_{{\cal Q}_{\epsilon_{j_0}} } \widetilde{w}_4= - \widetilde{g}/c $ and $
\|\widetilde{w}_4\|_1\leq 1.$  In fact, without
loss of generality, we assume that $ \{-\textbf{e}_i: ~i=1, \dots, q\}$
are arranged as the first $q$ columns  and  $ \{\textbf{e}_i: ~i=1,
\dots,  q\}$   as the second $q$ columns in
$\Gamma_{{\cal Q}_{\epsilon_{j_0}} }. $ For every $i=1, \dots , q,$ we set
$$ (\widetilde{w}_4)_i=\left\{\begin{array} {ll} (\widetilde{g})_i/c & \textrm{if} ~ (\widetilde{g})_i\geq 0,\\
 0 & \textrm{otherwise}, \end{array} \right.   ~  (\widetilde{w}_4)_{q+i} = \left\{\begin{array} {ll}  - (\widetilde{g})_i/c  & \textrm{if} ~(\widetilde{g})_i< 0, \\
 0 & \textrm{otherwise}. \end{array} \right.  $$
      All
remaining entries of $ \widetilde{w}_4  \in \mathbb{R}^{\widehat{n}}
$ are set to be zero.    By this choice, we see
that $\widetilde{w}_4 \geq 0$, $\Gamma_{{\cal Q}_{\epsilon_{j_0}} }
\widetilde{w}_4= - \widetilde{g} /c  $ and
\begin{eqnarray} \label{E1}  \widehat{\mathbf{e}}^T \widetilde{w}_4  &  = & \|\widetilde{w}_4\|_1  =   \|\widetilde{g}\|_1/c
 =  \|\widetilde{g}_\mathfrak{J} \|_1/c
    =    \|M_\mathfrak{J}^{-1} u^* \|_1/c \nonumber    =    \|M_\mathfrak{J}^{-1}(AA^T)^{-1}A \zeta\|_1/c \\
  & \leq &  \|M_\mathfrak{J}^{-1}(AA^T)^{-1}A\}_{\infty\to 1} \|\zeta\|_\infty/c
    \leq    \|\zeta\|_\infty =1.
\end{eqnarray}  By an argument similar to (\ref{6666}), we see that \begin{equation} \label{kiki} \phi^* (\widetilde{g})   \leq \|\widetilde{g} \|_1 \leq c. \end{equation}
  We also note that
$$ A^TM \Gamma_{{\cal Q}_{\epsilon_{j_0}}} \widetilde{w}_4 =A^TM(- \widetilde{g}/c )= A^T(-u^*/c) =-\zeta/c=-(\widetilde{w}_1-\widetilde{w}_2).   $$ Thus the vector $ (\widetilde{w}_1, \widetilde{w}_2, \widetilde{w}_3, \widetilde{w}_4) $ constructed   above satisfies the constraint of  (\ref{ddll}).
 Let $ \mathfrak{D}^{(2)} $ be given as  in Lemma 5.2.
 For the vector $(x,t,  \rho, \widetilde{w}_1, \widetilde{w}_2, \widetilde{w}_3, \widetilde{w}_4), $ by applying Lemma 2.1 with $ (M^1, M^2): = (\widehat{M}^1, \widehat{M}^2) $ where $\widehat{M}^1$ and $ \widehat{M}^2$ are given in  (\ref{M+b+}) and (\ref{M++}),
 we conclude that there is a point in $ \mathfrak{D}^{(2)} ,$ denoted by
$(\widehat{x}, \widehat{t}, \widehat{\rho},  \widehat{w}_1, \widehat{w}_2, \widehat{w}_3, \widehat{w}_4)), $ such that
\begin{equation} \label {hoff-4} \left\|
   \left[
      \begin{array}{c}
        x \\
        t \\
        \rho\\
        \widetilde{w}_1 \\
        \widetilde{w}_2 \\
        \widetilde{w}_3 \\
        \widetilde{w}_4 \\
      \end{array}\right] - \left[
       \begin{array}{c}
         \widehat{x} \\
        \widehat{t} \\
        \widehat{\rho} \\
         \widehat{w}_1 \\
         \widehat{w}_2 \\
         \widehat{w}_3 \\
           \widehat{w}_4 \\
       \end{array}
     \right]
 \right\|_2 \leq  \widehat{\gamma}  \left\| \left[\begin{array}{c}
(x-t)^+ \\
   (-x- t)^+ \\
   (\widetilde{\mathbf{e}}^Tt-\mu)^+\\
\left[(\Gamma_{{\cal Q}_{\epsilon_{j_0}} })^T[M^T(Ax-y)]-\rho \widehat{\mathbf{e}}  \right]^+ \\
    A^T M \Gamma_{{\cal Q}_{\epsilon_{j_0}} } \widetilde{w}_4+ \widetilde{w}_1-\widetilde{w}_2  \\
    ( \widetilde{w}_1+\widetilde{w}_2- \widetilde{w}_3 \widetilde{\mathbf{e}} )^+ \\
    ( \widehat{\mathbf{e}} ^T \widetilde{w}_4 -1)^+ \\
      \rho +  \mu \widetilde{w}_3 + ( y^TM \Gamma_{{\cal Q}_{\epsilon_{j_0}} })\widetilde{w}_4 \\
        (Z )^+
            \end{array}
             \right] \right\|_1, \end{equation}
    where  $$(Z)^+=
((-t)^+, ~(-\rho)^+, ~(-\widetilde{w}_1)^+,
          ~ (-\widetilde{w}_2)^+,
          ~ (-\widetilde{w}_3)^+, ~ (-\widetilde{w}_4)^+),
      $$ and  $\widehat{\gamma}  =\sigma_{\infty, 2} (\widehat{M}^1, \widehat{M}^2) $ is the Robinson's constant
    determined by  $(\widehat{M}^1, \widehat{M}^2)$ given in (\ref{M+b+}) and (\ref{M++}). The nonnegativity of $(t, \rho, \widetilde{w}_1, \widetilde{w}_2, \widetilde{w}_3, \widetilde{w}_4)$ implies that $   (Z)^+ =0.$
 The fact  $t=|x|$ implies that  $(x-t)^+ =
   (-x- t)^+    =0   $ and $\widetilde{\mathbf{e}}^T t =\|x\|_1.$  Since  $ (\widetilde{w}_1,
\widetilde{w}_2, \widetilde{w}_3, \widetilde{w}_4) $ is feasible to (\ref{ddll}), we
have      $$  ( \widehat{\mathbf{e}} ^T \widetilde{w}_4 -1)^+ =0,  ~ ( \widetilde{w}_1+\widetilde{w}_2- \widetilde{w}_3 \widetilde{\mathbf{e}})^+ =0, ~ A^T M \Gamma_{{\cal Q}_{\epsilon_{j_0}} } \widetilde{w}_4+ \widetilde{w}_1-\widetilde{w}_2 =0.$$
          Note that $\rho =\phi(M^T(Ax-y))$ implies (\ref{5544}), and hence   $ \left((\Gamma_{{\cal Q}_{\epsilon_{j_0}} })^T(M^T(Ax-y)) -\rho  \widehat{\mathbf{e}}  \right)^+ =0.   $
Thus it follows from   (\ref{hoff-4}) that
\begin{equation}  \label {error-A} \|x-\widehat{x} \|_2 \leq \widehat{\gamma} \left\{
  (\|x\|_1  - \mu)^+ +
\left|\rho+ \mu \widetilde{w}_3 +      y^T M \Gamma_{{\cal Q}_{\epsilon_{j_0}}
} \widetilde{w}_4\right|\right\}.
\end{equation}
From the definition of $\zeta  $  and $S,$   we have  $ x^T \zeta  = \|x_S\|_1+ x_{\overline{S}}^T
\zeta_{\overline{S}}   $  and $\|x_{\overline{S}}\|_1 =\sigma_k(x)_1 . $ Thus
$$ \left|  \|x\|_1 - x^T \zeta \right|= \left| \|x_{\overline{S}}\|_1 - x_{\overline{S}}^T
\zeta_{\overline{S}} \right|\leq     \| x_{\overline{S}}\|_1 + \|x_{\overline{S}} \|_1 \| \zeta_{\overline{S}}\| _\infty \leq 2\sigma_k(x)_1.$$
From (\ref{E1}), we see that $\|\widetilde{g}\|_1\leq c. $
 Note that  $\Gamma_{{\cal Q}_{\epsilon_{j_0}} }\widetilde{w}_4=-\widetilde{g} / c$ and $A^TM \Gamma_{{\cal Q}_{\epsilon_{j_0}}
} \widetilde{w}_4  = -(\widetilde{w}_1-\widetilde{w}_2) =-\zeta/c. $   By letting  $\psi = M^T(Ax-y) ,$   we have
 \begin{eqnarray} \label {error-B}   \left|\rho+ \mu \widetilde{w}_3 +      y^T M \Gamma_{{\cal Q}_{\epsilon_{j_0}}
} \widetilde{w}_4\right|
& = & \left| \rho +  \mu   \widetilde{w}_3      +  (x^T A^TM - \psi^T) \Gamma_{{\cal Q}_{\epsilon_{j_0}}
} \widetilde{w}_4 \right | \nonumber \\
& = &  \left|\rho +   \mu  \widetilde{w}_3   +  x^T A^TM \Gamma_{{\cal Q}_{\epsilon_{j_0}}
} \widetilde{w}_4  - \psi^T  \Gamma_{{\cal Q}_{\epsilon_{j_0}}
} \widetilde{w}_4 \right |  \nonumber \\
& \leq  &  \rho +  \frac{ \left|(\mu - \|x\|_1)   +   \|x\|_1     -   x^T \zeta   +  \psi^T  \widetilde{g} \right | } {c }     \nonumber \\
& \leq &  \rho +  \frac{ |\mu - \|x\|_1|   +  2\sigma_k (x)_1     +  \left|  \psi^T  \widetilde{g} \right | } {c }  \nonumber\\
& \leq &  \phi ( \psi)  +  \frac{ |\mu - \|x\|_1|   +  2\sigma_k (x)_1     +  \phi( \psi) \phi^*  ( \widetilde{g})   } {c }  \nonumber\\
& \leq  &   2\phi ( \psi)      +  \frac{|\mu - \|x\|_1|   +  2\sigma_k (x)_1 } {c   },
\end{eqnarray}
where the last two inequalities follow  from the fact $\rho= \phi ( \psi),$      $ \left|  \psi^T  \widetilde{g} \right | \leq  \phi( \psi) \phi^*  ( \widetilde{g}) $ and (\ref{kiki}).
Combining (\ref{error-A}) and  (\ref{error-B}) leads to
\begin{equation}\label{error-C} \|x-\widehat{x}\|_2 \leq
     \widehat{\gamma}
\left[ (\|x\|_1-\mu )^+  +   2\phi ( \psi)      +  \frac{|\mu - \|x\|_1|   +  2\sigma_k (x)_1 } {c   }  \right]. \end{equation}
  Let
$x^*$ and $ \overline{x} $ be the  projections of $x$
onto the compact convex sets $\Lambda^*$ and $\Lambda_{\epsilon_{j_0}},$ respectively, i.e.,
  $$ x^* = \Pi_ {\Lambda^*} (x)  \in \Lambda^*, ~  \overline{x}  =
\Pi_ {\Lambda_{\epsilon_{j_0}}}  (x) \in
\Lambda_{\epsilon_{j_0}}.  $$
By (\ref{DLT1}),   we have $ d^{\cal H} ( \Lambda^*, \Lambda_{\epsilon_{j_0}} )  \leq \delta.  $ Note that $\widehat{x} \in \Lambda^*_{\epsilon_{j_0}}  \subseteq  \Lambda_{\epsilon_{j_0}}  $ It follows from Lemma 4.1 that
$ \|x-x^*\|_2  \leq
    \delta + 2 \|x-\widehat{x}\|_2.  $
From this inequality and (\ref{error-C}), we conclude that
$$ \|x-x^*\|_2 \leq  \delta + 2  \widehat{\gamma} \left[
  (\|x\|_1-\mu )^+  + 2 \phi ( \psi)       +  \frac{|\mu - \|x\|_1|   +  2\sigma_k (x)_1 } {c   }    \right]. $$
Particularly, if $x $ obeys
$\|x\|_1 \leq \mu,$  we obtain
$$ \|x-x^*\|_2 \leq  \delta + 2 \widehat{\gamma}\left[
    2\phi ( \psi)       +  \frac{|\mu - \|x\|_1|   +  2\sigma_k (x)_1 } {c   }   \right], $$
as desired.
 ~~ $\Box $

\vskip 0.07in

When the parameter $\mu $ is  large,  the optimal solution $x^*$ of (\ref{GLASSO}) might be taken in the interior of the feasible set, i.e., $\|x^*\|_1< \mu. $    When $\mu $ is   small, the optimal solution of  (\ref{GLASSO}) usually attains at the boundary of its feasible set, i.e., $\|x^*\|_1=\mu.$ Thus, in stability analysis of LASSO, we are particularly interested in the gap between $x^*$ and those vectors satisfying   $ \|x\|_1 <  \mu$ or $ \|x\|_1 = \mu.$
The following  result    follows immediately from Theorem 5.3.

\vskip 0.07in

 \textbf{ Corollary 5.4.} \emph{Let $\delta >0$ be any fixed   small constant, and let $\mathcal{Q}_{\epsilon_{j_0}}$ be the fixed polytope represented as (\ref{QQQ}) such that (\ref{DLT1}) is achieved. Let the   data $(M, A, y,
\mu)$  in (\ref{GLASSO}) be given, where $\mu >0 , $
  $ A \in \mathbb{R}^{m\times n}$  ($m<n$)  and $M \in R^{m\times q} $ ($m\leq q$)  with $rank (A)=rank (M)=m.$
Suppose that one of the conditions listed in Corollary 3.3 is satisfied. Then the following statements hold:  }

(i)  \emph{For any $x\in
\mathbb{R}^n$ with $ \|x\|_1 <  \mu , $ there is an
optimal solution $x^*$ of (\ref{GLASSO})
     approximating  $x$ with  error (\ref{1616}).  }

(ii)   \emph{For any $x\in
\mathbb{R}^n$ with  $ \|x\|_1 =\mu , $ there is an
optimal solution $x^*$ of (\ref{GLASSO})
     approximating  $x$ with error
\begin{equation} \label{FFNN}  \|x-x^*\|_2 \leq  \delta + 4    \widehat{\gamma}
\left[  \phi (M^T (Ax-y))      +  \frac{ \sigma_k (x)_1 } {c   }  \right], \end{equation}
   where   the constants $c $ and  $ \widehat{\gamma} $ are given as in Theorem 5.3.
}

 \vskip 0.07in

  The   stability results for the special cases   $M=I$ or $M=A$ can be obtained immediately from the above result. The statement of such  a result  is omitted here.

\vskip 0.15in

\section{Conclusions} We have shown that the general Dantzig selector and LASSO problems are   stable in sparse data recovery  under the so-called weak range space property of a transposed design matrix. These optimization problems  are general enough to include many important special cases.  Our   stability analysis for these problems is carried out   differently from those in the literature in terms of the analytic method, mild assumption and the way of expression of stability coefficients. The classic Hoffman's Lemma and a polytope approximation technique of convex bodies are employed as a new deterministic analytic method for the development of a stability theory for Dantzig selector and LASSO problems.    The stability coefficients are measured by   Robinson's constants depending on the problem data. The assumption made in this paper is a constant-free  matrix condition, which is  naturally originated from the fundamental optimality condition of convex optimization. It turns out to be a    mild sufficient condition for Dantzig selector and LASSO  problems to be stable in sparse data recovery.  We have shown that this assumption  is   also necessary for the standard Dantzig selector   to be  stable. Many known matrix  conditions in compressed sensing such as  RIP, NSP and others imply our assumption.


\baselineskip 0.18in
\begin{thebibliography}{999}


\bibitem{AS14} J. Andersson and J.O. Str\"omberg, On the theorem of
uniform recovery of structured random matrices, \emph{IEEE Trans. Inform. Theory},
60 (2014), pp. 1700--1710. \\[-5.5mm]


\bibitem{BDMS13}  A. Bandeira, E. Dobriban, D. Mixon and W. Sawin, Certifying the restricted
isometry property is hard, \emph{IEEE Trans. Inform. Theory},  59 (2013),
pp. 3448--3450.    \\[-5.5mm]

\bibitem{B14} A. Barvinok, Thrifty approximations of convex bodies by polytopes, \emph{Int. Math. Res. Notices},  16 (2014),  pp. 4341--4356.  \\[-5.5mm]



\bibitem{BRT09} P. Bickel, Y. Ritov, and A. Tsybakov, Simultaneous analysis of lasso and dantzig selector,
\emph{Ann. Statist.}, 37 (2009), pp. 1705-1732.  \\[-5.5mm]


\bibitem{BI75} E.M. Bronshtein and L. D. Ivanov. The approximation of convex sets by polyhedra,
   \emph{Siberian Math. J.}, 16 (1975),     pp. 852--853.   \\[-5.5mm]

   \bibitem{BDE09} A.M. Bruckstein, D. Donoho and M. Elad,
From sparse solutions of systems of equations to sparse modeling of
signals and images, \emph{SIAM Rev.}, 51 (2009), pp. 34--81.   \\[-5.5mm]



  \bibitem{BG11} P. B\"uhlmann and S. van de Geer,   \emph{Statistics for High-Dimensional Data:  Methods, Theory and Applications}, Springer, Berlin, 2011.   \\[-5.5mm]




\bibitem{CCW15} J. Cahill, X. Chen and R. Wang,  The gap between the
null space property and the restricted isometry property, \emph{Linear Algebra Appl.},  501 (2016),   pp. 363--375.   \\[-5.5mm]

\bibitem{CWX13} T. Cai and A. Zhang, Sharp RIP bound for sparse
signal and low-rank matrix recovery,  \emph{Appl.  Comput. Harmon.
Anal.}, 35 (2013),  pp. 74--93.   \\[-5.5mm]

\bibitem{CZ14} T. Cai and A. Zhang, Sparse representation of a
polytope and recovery of sparse signals and low-rank matrices,
\emph{IEEE Trans. Inform. Theory}, 60 (2014), pp. 122--132. \\[-5.5mm]

\bibitem{CXZ09} T. Cai, G. Xu, J. Zhang, On recovery of sparse signals via $\ell_1$ minimization,  \emph{IEEE Trans. Inform. Theory}, 55 (2009), pp.3388--3397.   \\[-5.5mm]



\bibitem{CT07} E. Candes and T. Tao, The Dantzig selector: statistical estimation when $p$ is much larger
than $n,$ \emph{Ann. Statist.},  35 (2007),  pp. 2313--2351.  \\[-5.5mm]

\bibitem{CT07b} E. Candes and T. Tao, Rejoinder: the Dantzig selector: statistical
estimation when p is much larger than n, \emph{Ann. Statist.},
  35 (2007), pp. 2392--2404.  \\[-5.5mm]

\bibitem{CRT06b} E. Cand$\grave{\textrm{e}}$s, J. Romberg  and T.
Tao,  Stable signal recovery from incomplete and inacurate
measurements, \emph{Comm. Pure Appl. Math.}, 59 (2006), pp.
1207--1223.   \\[-5.5mm]

\bibitem{CT05} E. Cand$\grave{\textrm{e}}$s and T. Tao, Decoding by linear
programming, \emph{IEEE Trans. Inform. Theory}, 51 (2005), pp.
4203--4215.   \\[-5.5mm]

\bibitem{CWB08} E. Cand$\grave{\textrm{e}}$s, M. Wakin and S.
Boyd, Enhancing sparsity by reweighted $\ell_1$ minimization,
\emph{J. Fourier Anal. Appl.}, 14 (2008), pp. 877--905.   \\[-5.5mm]


\bibitem{CDS98} S. Chen, D. Donoho and M. Saunders, Atomic
decomposition by basis pursuit, \emph{SIAM J. Sci. Comput.}, 20 (1998),
pp. 33--61.   \\[-5.5mm]

\bibitem{CDD09} A. Cohen, W. Dahmen and R. Devore, Compressed
sensing and best $k$-term aproximation, \emph{J. Amer. Math. Soc.},
22 (2009), pp. 211--231.   \\[-5.5mm]





\bibitem{C12} Y. De Castro, A remark on the LASSO and the Dantzig selector, \emph{Stat. Prob. Lett.}, 83 (2013), pp. 304--314.   \\[-5.5mm]

\bibitem{D06} D. L. Donoho, Compressed sensing, \emph{IEEE Trans. Inform.
Theory}, 52 (2006), pp. 1289--1306.    \\[-5.5mm]

\bibitem{DET06} D. Donoho, M. Elad and V. Temlyahov, Stable recovery
of sparse overcomplete representations in the presenxe of noise,
\emph{IEEE Trans. Inform. Theory,} 52 (2006), pp. 6--18.    \\[-5.5mm]

\bibitem{D74} R. Dudley, Matric entropy of some classes of sets with
differentiable bounaries, \emph{J. Approx. Theory, } 10 (1974),
227-236; Correction, \emph{J. Approx. Theory}, 26 (1979), pp.192--193.    \\[-5.5mm]

\bibitem{EHT07} B. Efron, T. Hastie, and R. Tibshirani, Discussion: the Dantzig
selector: statistical estimation when p is much larger than n, \emph{Ann.  Statist.},  35 (2007), pp. 2358--2364.   \\[-5.5mm]

\bibitem{E10} M. Elad, \emph{Sparse and Redundant Representations: From Theory to Applications in Signal
and Image Processing}, Springer, New York, 2010.    \\[-5.5mm]

\bibitem{EK12} Y. Eldar and G. Kutyniok, \emph{Compressed Sensing: Theory
and Applications}, Cambridge University Press, 2012.    \\[-5.5mm]

\bibitem{FS07} M. Friedlander and M. Saunders , Discussion: the Dantzig selector:
statistical estimation when p is much larger than n, \emph{Ann.
Statist.}, 35 (2007), pp. 2385--2391.   \\[-5.5mm]


\bibitem{F14} S. Foucart,  Stability and robustness of
$\ell_1$-minimization with Weibull matrices and redandant
distionaries, \emph{Linear Algebra Appl.}, 441
(2014), pp. 4--21.   \\[-5.5mm]

\bibitem{FR13} S. Foucart and H. Rauhut, \emph{A Mathematical Introduction to Compressive Sensing},
  Springer, NY, 2013.   \\[-5.5mm]

\bibitem{F04} J. Fuchs, On sparse representations in arbitrary
redundant bases, \emph{IEEE Trans. Inform. Theory,} 50 (2004),
pp. 1341--1344.   \\[-5.5mm]

\bibitem{GSH11} M. Grasmair, O. Sherzer and M. Haltmeier, Necessary and sufficient conditions for
linear convergence of l1-regularization, \emph{Comm. Pure Appl.
Math.}, 64 (2011), pp. 161--182.    \\[-5.5mm]

\bibitem{GB09} S.A. van de Geer and P. B\"uhlmann, On the condition used to prove oracle results for the Lasso, \emph{Electron. J. Stat.}, 3 (2009), pp. 1360--1392.   \\[-5.5mm]

  \bibitem{HTW15} T. Hastie, R. Tibshirani and M. Wainwright,  \emph{Statistical Learning with Sparsity: The Lasso and Generalizations}, Chapman \&  Hall/CRC,  2015.   \\[-5.5mm]

\bibitem{H52} A. J. Hoffman, On the approximation solution of
systems of linear inequalities, \emph{J. Res. Nat. Bur. Standards},
49 (1952), pp. 263--265.    \\[-5.5mm]


\bibitem{JRL09} G. M. James, P. Radchenko, and J. Lv,  Dasso: connections between the dantzig selector
and lasso, \emph{J. Roy. Statist. Soc. Ser. B}, 71 (2009),
pp. 127--142.    \\[-5.5mm]

\bibitem{J15} O. James, Stabibility analysis of LASSO and Dantzig selector via constrained minimal sigular value of Gaussian sensing matrices, International Conference on Current Trends in Advanced Computing (ICCTAC-2015) ICCTAC 2015(2): 1-5, May 2015.     \\[-5.5mm]

\bibitem{JN11} A. Juditsky and A. Nemirovski, Accuracy guarantees for $\ell_1$-recovery, \emph{IEEE Trans.  Inform. Theory},   57 (2011), pp. 7818--7839.    \\[-5.5mm]

\bibitem{KZL15} Y. Kong, Z. Zheng, and J. Lv, The constrained Dantzig selector with enhanced consistency, \emph{J. Machine Learning Res.},  17 (2016) pp. 1--22.    \\[-5.5mm]

\bibitem {KT51} H. W. Kuhn,  and  A. W. Tucker, \emph{Nonlinear programming},  Proceedings of 2nd Berkeley Symposium,
 University of California Press, Berkeley, pp. 481--492, 1951.     \\[-5.5mm]



\bibitem{MS87} O.L. Mangasarian and T. Shiau, Lipschitz continuity of solutions of linear inequalities, programs and complementarity problems, \emph{SIAM J. Control \& Optim.}, 25(1987), pp. 583--595.    \\[-5.5mm]

\bibitem{MR15} R Mazumder and P Radchenko, The discrete Dantzig selector: Estimating sparse linear models via mixed integer linear optimization, arXiv:1508.01922v3, 2017.    \\[-5.5mm]

\bibitem{MRY07}  N. Meinshausen, G. Rocha, and B. Yu, Discussion: a tale of three
cousins: lasso, L2boosting and Dantzig, \emph{Ann. Statist.},
35 (2007), pp. 2373--2384.   \\[-5.5mm]


\bibitem{MB06} N. Meinshausen and P. B\"uhlmann, Stability selection, \emph{J. Royal Stat. Soc. Ser. B}, 72 (2010), pp. 417--473.     \\[-5.5mm]




\bibitem{P89} G. Pisier,  \emph{The volume of convex bodies and Banach space geometry}, Cambridge Tracts in Mathematics 94, Cambridge University Press, 1989.   \\[-5.5mm]

\bibitem{P07} M. Plumbley, On polar polytopes and the recovery of sparse representations,
\emph{IEEE Trans. Inform. Theory}, 53 (2007), pp. 3188--3195.    \\[-5.5mm]


\bibitem{R73} S. M. Robinson, Bounds for error in the solution set
of a perturbed linear program, \emph{Linear Algebra  Appl.}, 6
(1973), pp. 69--81.    \\[-5.5mm]

\bibitem {R70} R.T. Rockafellar, \emph{Convex Analysis}, Princeton University Press, Princeton, New Jersey, 1970.    \\[-5.5mm]

 \bibitem{RZ13} M. Rudelson and S. Zhou, Reconstruction from anisotropic random measurements, \emph{IEEE Trans. Inform. Theory}, 59 (2013), pp. 3434--3447.    \\[-5.5mm]

\bibitem{AR09} M. Salman Asif and J. Romberg, On the LASSO and Dantzig selector equivalence,   in \emph{IEEE 44th Annual Conference on Information Sciences and Systems}, 2010,  pp. 1--6.   \\[-5.5mm]



\bibitem{T96} R. Tibshirani,  Regression shrinkage and selection via the lasso. \emph{J. Roy. Statist. Soc. Ser. B}, 58 (1996), pp. 267--288.   \\[-5.5mm]


\bibitem{TP14} A. Tillmann and M. Pfetsch, The computational complexity of the restricted isometry property, the nullspace property, and related concepts in compressed sensing,
\emph{IEEE Trans. Inform. Theory},  60 (2014), pp. 1248--1259.   \\[-5.5mm]


\bibitem{T04} J.A. Tropp,   Greed is good: Algorithmic results for sparse approximation,
\emph{IEEE Trans. Inform. Theory}, 50 (2004), pp. 2231--2242.   \\[-5.5mm]

\bibitem {K39} W. Karush,  Minima of functions of several variables with inequalities as side constraints,   MSc Dissertation,  Univ. of Chicago, Chicago, 1939.   \\[-5.5mm]


\bibitem{ZY06} P. Zhao and B. Yu,  On model selection consistency of Lasso, \emph{J. Machine Learning Res.}, 7 (2006), pp. 2541--2567.   \\[-5.5mm]

\bibitem{Z13} Y.-B. Zhao, RSP-based analysis for sparsest
and least $\ell_1$-norm  solutions to
 underdetermined  linear systems, \emph{IEEE Trans. Signal
 Process.}, 61 (2013), no. 22,  pp. 5777--5788.    \\[-5.5mm]



\bibitem{ZL12} Y.-B. Zhao and D. Li,   Reweighted $\ell_1$-minimization for sparse solutions to
underdetermined linear systems, \emph{SIAM J. Optim.}, 22 (2012), pp. 893--912.    \\[-5.5mm]


 \bibitem{ZK15} Y.-B. Zhao and M. Kocvara, A new computational  method for the sparsest solutions to
  systems of linear equations, \emph{SIAM J. Optim.},    25 (2015),  pp. 1110--1134.   \\[-5.5mm]


\bibitem{ZX14}   Y.-B. Zhao and C. Xu,  1-bit compressive
 sensing: Reformulation and RRSP-based sign recovery theory,   \emph{Sci. China Math.}, 59 (2016),   pp. 2049--2074.   \\[-5.5mm]

\end{thebibliography}
\end{document}